\documentclass[11pt]{article}
\usepackage{amsfonts}
\long\def\ig#1{\relax}
\ig{Thanks to Roberto Minio for this def'n.  Compare the def'n of
\comment in AMSTeX.}

\newcount \coefa
\newcount \coefb
\newcount \coefc
\newcount\tempcounta
\newcount\tempcountb
\newcount\tempcountc
\newcount\tempcountd
\newcount\xext
\newcount\yext
\newcount\xoff
\newcount\yoff
\newcount\gap%
\newcount\arrowtypea
\newcount\arrowtypeb
\newcount\arrowtypec
\newcount\arrowtyped
\newcount\arrowtypee
\newcount\height
\newcount\width
\newcount\xpos
\newcount\ypos
\newcount\run
\newcount\rise
\newcount\arrowlength
\newcount\halflength
\newcount\arrowtype
\newdimen\tempdimen
\newdimen\xlen
\newdimen\ylen
\newsavebox{\tempboxa}%
\newsavebox{\tempboxb}%
\newsavebox{\tempboxc}%

\makeatletter
\def\settypes(#1,#2,#3){\arrowtypea#1 \arrowtypeb#2 \arrowtypec#3}
\def\settoheight#1#2{\setbox\@tempboxa\hbox{#2}#1\ht\@tempboxa\relax}%
\def\settodepth#1#2{\setbox\@tempboxa\hbox{#2}#1\dp\@tempboxa\relax}%
\def\settokens[#1`#2`#3`#4]{%
     \def\tokena{#1}\def\tokenb{#2}\def\tokenc{#3}\def\tokend{#4}}
\def\setsqparms[#1`#2`#3`#4;#5`#6]{%
\arrowtypea #1
\arrowtypeb #2
\arrowtypec #3
\arrowtyped #4
\width #5
\height #6
}
\def\setpos(#1,#2){\xpos=#1 \ypos#2}

\def\bfig{\begin{picture}(\xext,\yext)(\xoff,\yoff)}
\def\efig{\end{picture}}

\def\putbox(#1,#2)#3{\put(#1,#2){\makebox(0,0){$#3$}}}

\def\settriparms[#1`#2`#3;#4]{\settripairparms[#1`#2`#3`1`1;#4]}%

\def\settripairparms[#1`#2`#3`#4`#5;#6]{%
\arrowtypea #1
\arrowtypeb #2
\arrowtypec #3
\arrowtyped #4
\arrowtypee #5
\width #6
\height #6
}

\def\resetparms{\settripairparms[1`1`1`1`1;500]\width 500}

\resetparms

\def\mvector(#1,#2)#3{
\put(0,0){\vector(#1,#2){#3}}%
\put(0,0){\vector(#1,#2){30}}%
}
\def\evector(#1,#2)#3{{
\arrowlength #3
\put(0,0){\vector(#1,#2){\arrowlength}}%
\advance \arrowlength by-30
\put(0,0){\vector(#1,#2){\arrowlength}}%
}}

\def\horsize#1#2{%
\settowidth{\tempdimen}{$#2$}%
#1=\tempdimen
\divide #1 by\unitlength
}

\def\vertsize#1#2{%
\settoheight{\tempdimen}{$#2$}%
#1=\tempdimen
\settodepth{\tempdimen}{$#2$}%
\advance #1 by\tempdimen
\divide #1 by\unitlength
}

\def\vertadjust[#1`#2`#3]{%
\vertsize{\tempcounta}{#1}%
\vertsize{\tempcountb}{#2}%
\ifnum \tempcounta<\tempcountb \tempcounta=\tempcountb \fi
\divide\tempcounta by2
\vertsize{\tempcountb}{#3}%
\ifnum \tempcountb>0 \advance \tempcountb by20 \fi
\ifnum \tempcounta<\tempcountb \tempcounta=\tempcountb \fi
}

\def\horadjust[#1`#2`#3]{%
\horsize{\tempcounta}{#1}%
\horsize{\tempcountb}{#2}%
\ifnum \tempcounta<\tempcountb \tempcounta=\tempcountb \fi
\divide\tempcounta by20
\horsize{\tempcountb}{#3}%
\ifnum \tempcountb>0 \advance \tempcountb by60 \fi
\ifnum \tempcounta<\tempcountb \tempcounta=\tempcountb \fi
}

\ig{ In this procedure, #1 is the paramater that sticks out all the way,
#2 sticks out the least and #3 is a label sticking out half way.  #4 is
the amount of the offset.}

\def\sladjust[#1`#2`#3]#4{%
\tempcountc=#4
\horsize{\tempcounta}{#1}%
\divide \tempcounta by2
\horsize{\tempcountb}{#2}%
\divide \tempcountb by2
\advance \tempcountb by-\tempcountc
\ifnum \tempcounta<\tempcountb \tempcounta=\tempcountb\fi
\divide \tempcountc by2
\horsize{\tempcountb}{#3}%
\advance \tempcountb by-\tempcountc
\ifnum \tempcountb>0 \advance \tempcountb by80\fi
\ifnum \tempcounta<\tempcountb \tempcounta=\tempcountb\fi
\advance\tempcounta by20
}

\def\putvector(#1,#2)(#3,#4)#5#6{{%
\xpos=#1
\ypos=#2
\run=#3
\rise=#4
\arrowlength=#5
\arrowtype=#6
\ifnum \arrowtype<0
    \ifnum \run=0
        \advance \ypos by-\arrowlength
    \else
        \tempcounta \arrowlength
        \multiply \tempcounta by\rise
        \divide \tempcounta by\run
        \ifnum\run>0
            \advance \xpos by\arrowlength
            \advance \ypos by\tempcounta
        \else
            \advance \xpos by-\arrowlength
            \advance \ypos by-\tempcounta
        \fi
    \fi
    \multiply \arrowtype by-1
    \multiply \rise by-1
    \multiply \run by-1
\fi
\ifnum \arrowtype=1
    \put(\xpos,\ypos){\vector(\run,\rise){\arrowlength}}%
\else\ifnum \arrowtype=2
    \put(\xpos,\ypos){\mvector(\run,\rise)\arrowlength}%
\else\ifnum\arrowtype=3
    \put(\xpos,\ypos){\evector(\run,\rise){\arrowlength}}%
\fi\fi\fi
}}

\def\putsplitvector(#1,#2)#3#4{
\xpos #1
\ypos #2
\arrowtype #4
\halflength #3
\arrowlength #3
\gap 140
\advance \halflength by-\gap
\divide \halflength by2
\ifnum \arrowtype=1
    \put(\xpos,\ypos){\line(0,-1){\halflength}}%
    \advance\ypos by-\halflength
    \advance\ypos by-\gap
    \put(\xpos,\ypos){\vector(0,-1){\halflength}}%
\else\ifnum \arrowtype=2
    \put(\xpos,\ypos){\line(0,-1)\halflength}%
    \put(\xpos,\ypos){\vector(0,-1)3}%
    \advance\ypos by-\halflength
    \advance\ypos by-\gap
    \put(\xpos,\ypos){\vector(0,-1){\halflength}}%
\else\ifnum\arrowtype=3
    \put(\xpos,\ypos){\line(0,-1)\halflength}%
    \advance\ypos by-\halflength
    \advance\ypos by-\gap
    \put(\xpos,\ypos){\evector(0,-1){\halflength}}%
\else\ifnum \arrowtype=-1
    \advance \ypos by-\arrowlength
    \put(\xpos,\ypos){\line(0,1){\halflength}}%
    \advance\ypos by\halflength
    \advance\ypos by\gap
    \put(\xpos,\ypos){\vector(0,1){\halflength}}%
\else\ifnum \arrowtype=-2
    \advance \ypos by-\arrowlength
    \put(\xpos,\ypos){\line(0,1)\halflength}%
    \put(\xpos,\ypos){\vector(0,1)3}%
    \advance\ypos by\halflength
    \advance\ypos by\gap
    \put(\xpos,\ypos){\vector(0,1){\halflength}}%
\else\ifnum\arrowtype=-3
    \advance \ypos by-\arrowlength
    \put(\xpos,\ypos){\line(0,1)\halflength}%
    \advance\ypos by\halflength
    \advance\ypos by\gap
    \put(\xpos,\ypos){\evector(0,1){\halflength}}%
\fi\fi\fi\fi\fi\fi
}

\def\putmorphism(#1)(#2,#3)[#4`#5`#6]#7#8#9{{%
\run #2
\rise #3
\ifnum\rise=0
  \puthmorphism(#1)[#4`#5`#6]{#7}{#8}{#9}%
\else\ifnum\run=0
  \putvmorphism(#1)[#4`#5`#6]{#7}{#8}{#9}%
\else
\setpos(#1)%
\arrowlength #7
\arrowtype #8
\ifnum\run=0
\else\ifnum\rise=0
\else
\ifnum\run>0
    \coefa=1
\else
   \coefa=-1
\fi
\ifnum\arrowtype>0
   \coefb=0
   \coefc=-1
\else
   \coefb=\coefa
   \coefc=1
   \arrowtype=-\arrowtype
\fi
\width=2
\multiply \width by\run
\divide \width by\rise
\ifnum \width<0  \width=-\width\fi
\advance\width by60
\if l#9 \width=-\width\fi
\putbox(\xpos,\ypos){#4}
{\multiply \coefa by\arrowlength
\advance\xpos by\coefa
\multiply \coefa by\rise
\divide \coefa by\run
\advance \ypos by\coefa
\putbox(\xpos,\ypos){#5} }%
{\multiply \coefa by\arrowlength
\divide \coefa by2
\advance \xpos by\coefa
\advance \xpos by\width
\multiply \coefa by\rise
\divide \coefa by\run
\advance \ypos by\coefa
\if l#9%
   \put(\xpos,\ypos){\makebox(0,0)[r]{$#6$}}%
\else\if r#9%
   \put(\xpos,\ypos){\makebox(0,0)[l]{$#6$}}%
\fi\fi }%
{\multiply \rise by-\coefc
\multiply \run by-\coefc
\multiply \coefb by\arrowlength
\advance \xpos by\coefb
\multiply \coefb by\rise
\divide \coefb by\run
\advance \ypos by\coefb
\multiply \coefc by70
\advance \ypos by\coefc
\multiply \coefc by\run
\divide \coefc by\rise
\advance \xpos by\coefc
\multiply \coefa by140
\multiply \coefa by\run
\divide \coefa by\rise
\advance \arrowlength by\coefa
\ifnum \arrowtype=1
   \put(\xpos,\ypos){\vector(\run,\rise){\arrowlength}}%
\else\ifnum\arrowtype=2
   \put(\xpos,\ypos){\mvector(\run,\rise){\arrowlength}}%
\else\ifnum\arrowtype=3
   \put(\xpos,\ypos){\evector(\run,\rise){\arrowlength}}%
\fi\fi\fi}\fi\fi\fi\fi}}

\def\puthmorphism(#1,#2)[#3`#4`#5]#6#7#8{{%
\xpos #1
\ypos #2
\width #6
\arrowlength #6
\putbox(\xpos,\ypos){#3\vphantom{#4}}%
{\advance \xpos by\arrowlength
\putbox(\xpos,\ypos){\vphantom{#3}#4}}%
\horsize{\tempcounta}{#3}%
\horsize{\tempcountb}{#4}%
\divide \tempcounta by2
\divide \tempcountb by2
\advance \tempcounta by30
\advance \tempcountb by30
\advance \xpos by\tempcounta
\advance \arrowlength by-\tempcounta
\advance \arrowlength by-\tempcountb
\putvector(\xpos,\ypos)(1,0){\arrowlength}{#7}%
\divide \arrowlength by2
\advance \xpos by\arrowlength
\vertsize{\tempcounta}{#5}%
\divide\tempcounta by2
\advance \tempcounta by20
\if a#8 %
   \advance \ypos by\tempcounta
   \putbox(\xpos,\ypos){#5}%
\else
   \advance \ypos by-\tempcounta
   \putbox(\xpos,\ypos){#5}%
\fi}}

\def\putvmorphism(#1,#2)[#3`#4`#5]#6#7#8{{%
\xpos #1
\ypos #2
\arrowlength #6
\arrowtype #7
\settowidth{\xlen}{$#5$}%
\putbox(\xpos,\ypos){#3}%
{\advance \ypos by-\arrowlength
\putbox(\xpos,\ypos){#4}}%
{\advance\arrowlength by-140
\advance \ypos by-70
\ifdim\xlen>0pt
   \if m#8%
      \putsplitvector(\xpos,\ypos){\arrowlength}{\arrowtype}%
   \else
      \putvector(\xpos,\ypos)(0,-1){\arrowlength}{\arrowtype}%
   \fi
\else
   \putvector(\xpos,\ypos)(0,-1){\arrowlength}{\arrowtype}%
\fi}%
\ifdim\xlen>0pt
   \divide \arrowlength by2
   \advance\ypos by-\arrowlength
   \if l#8%
      \advance \xpos by-40
      \put(\xpos,\ypos){\makebox(0,0)[r]{$#5$}}%
   \else\if r#8%
      \advance \xpos by40
      \put(\xpos,\ypos){\makebox(0,0)[l]{$#5$}}%
   \else
      \putbox(\xpos,\ypos){#5}%
   \fi\fi
\fi
}}

\def\topadjust[#1`#2`#3]{%
\yoff=10
\vertadjust[#1`#2`{#3}]%
\advance \yext by\tempcounta
\advance \yext by 10
}
\def\botadjust[#1`#2`#3]{%
\vertadjust[#1`#2`{#3}]%
\advance \yext by\tempcounta
\advance \yoff by-\tempcounta
}
\def\leftadjust[#1`#2`#3]{%
\xoff=0
\horadjust[#1`#2`{#3}]%
\advance \xext by\tempcounta
\advance \xoff by-\tempcounta
}
\def\rightadjust[#1`#2`#3]{%
\horadjust[#1`#2`{#3}]%
\advance \xext by\tempcounta
}
\def\rightsladjust[#1`#2`#3]{%
\sladjust[#1`#2`{#3}]{\width}%
\advance \xext by\tempcounta
}
\def\leftsladjust[#1`#2`#3]{%
\xoff=0
\sladjust[#1`#2`{#3}]{\width}%
\advance \xext by\tempcounta
\advance \xoff by-\tempcounta
}
\def\adjust[#1`#2;#3`#4;#5`#6;#7`#8]{%
\topadjust[#1``{#2}]
\leftadjust[#3``{#4}]
\rightadjust[#5``{#6}]
\botadjust[#7``{#8}]}

\def\putsquarep<#1>(#2)[#3;#4`#5`#6`#7]{{%
\setsqparms[#1]%
\setpos(#2)%
\settokens[#3]%
\puthmorphism(\xpos,\ypos)[\tokenc`\tokend`{#7}]{\width}{\arrowtyped}b%
\advance\ypos by \height
\puthmorphism(\xpos,\ypos)[\tokena`\tokenb`{#4}]{\width}{\arrowtypea}a%
\putvmorphism(\xpos,\ypos)[``{#5}]{\height}{\arrowtypeb}l%
\advance\xpos by \width
\putvmorphism(\xpos,\ypos)[``{#6}]{\height}{\arrowtypec}r%
}}

\def\putsquare{\@ifnextchar <{\putsquarep}{\putsquarep%
   <\arrowtypea`\arrowtypeb`\arrowtypec`\arrowtyped;\width`\height>}}
\def\square{\@ifnextchar< {\squarep}{\squarep
   <\arrowtypea`\arrowtypeb`\arrowtypec`\arrowtyped;\width`\height>}}
\def\squarep<#1>[#2`#3`#4`#5;#6`#7`#8`#9]{{
\setsqparms[#1]
\xext=\width                                          
\yext=\height                                         
\topadjust[#2`#3`{#6}]
\botadjust[#4`#5`{#9}]
\leftadjust[#2`#4`{#7}]
\rightadjust[#3`#5`{#8}]
\begin{picture}(\xext,\yext)(\xoff,\yoff)
\putsquarep<\arrowtypea`\arrowtypeb`\arrowtypec`\arrowtyped;\width`\height>%
(0,0)[#2`#3`#4`#5;#6`#7`#8`{#9}]%
\end{picture}%
}}

\def\putptrianglep<#1>(#2,#3)[#4`#5`#6;#7`#8`#9]{{%
\settriparms[#1]%
\xpos=#2 \ypos=#3
\advance\ypos by \height
\puthmorphism(\xpos,\ypos)[#4`#5`{#7}]{\height}{\arrowtypea}a%
\putvmorphism(\xpos,\ypos)[`#6`{#8}]{\height}{\arrowtypeb}l%
\advance\xpos by\height
\putmorphism(\xpos,\ypos)(-1,-1)[``{#9}]{\height}{\arrowtypec}r%
}}

\def\putptriangle{\@ifnextchar <{\putptrianglep}{\putptrianglep
   <\arrowtypea`\arrowtypeb`\arrowtypec;\height>}}
\def\ptriangle{\@ifnextchar <{\ptrianglep}{\ptrianglep
   <\arrowtypea`\arrowtypeb`\arrowtypec;\height>}}

\def\ptrianglep<#1>[#2`#3`#4;#5`#6`#7]{{
\settriparms[#1]%
\width=\height                         
\xext=\width                           
\yext=\width                           
\topadjust[#2`#3`{#5}]
\botadjust[#3``]
\leftadjust[#2`#4`{#6}]
\rightsladjust[#3`#4`{#7}]
\begin{picture}(\xext,\yext)(\xoff,\yoff)
\putptrianglep<\arrowtypea`\arrowtypeb`\arrowtypec;\height>%
(0,0)[#2`#3`#4;#5`#6`{#7}]%
\end{picture}%
}}

\def\putqtrianglep<#1>(#2,#3)[#4`#5`#6;#7`#8`#9]{{%
\settriparms[#1]%
\xpos=#2 \ypos=#3
\advance\ypos by\height
\puthmorphism(\xpos,\ypos)[#4`#5`{#7}]{\height}{\arrowtypea}a%
\putmorphism(\xpos,\ypos)(1,-1)[``{#8}]{\height}{\arrowtypeb}l%
\advance\xpos by\height
\putvmorphism(\xpos,\ypos)[`#6`{#9}]{\height}{\arrowtypec}r%
}}

\def\putqtriangle{\@ifnextchar <{\putqtrianglep}{\putqtrianglep
   <\arrowtypea`\arrowtypeb`\arrowtypec;\height>}}
\def\qtriangle{\@ifnextchar <{\qtrianglep}{\qtrianglep
   <\arrowtypea`\arrowtypeb`\arrowtypec;\height>}}

\def\qtrianglep<#1>[#2`#3`#4;#5`#6`#7]{{
\settriparms[#1]
\width=\height                         
\xext=\width                           
\yext=\height                          
\topadjust[#2`#3`{#5}]
\botadjust[#4``]
\leftsladjust[#2`#4`{#6}]
\rightadjust[#3`#4`{#7}]
\begin{picture}(\xext,\yext)(\xoff,\yoff)
\putqtrianglep<\arrowtypea`\arrowtypeb`\arrowtypec;\height>%
(0,0)[#2`#3`#4;#5`#6`{#7}]%
\end{picture}%
}}

\def\putdtrianglep<#1>(#2,#3)[#4`#5`#6;#7`#8`#9]{{%
\settriparms[#1]%
\xpos=#2 \ypos=#3
\puthmorphism(\xpos,\ypos)[#5`#6`{#9}]{\height}{\arrowtypec}b%
\advance\xpos by \height \advance\ypos by\height
\putmorphism(\xpos,\ypos)(-1,-1)[``{#7}]{\height}{\arrowtypea}l%
\putvmorphism(\xpos,\ypos)[#4``{#8}]{\height}{\arrowtypeb}r%
}}

\def\putdtriangle{\@ifnextchar <{\putdtrianglep}{\putdtrianglep
   <\arrowtypea`\arrowtypeb`\arrowtypec;\height>}}
\def\dtriangle{\@ifnextchar <{\dtrianglep}{\dtrianglep
   <\arrowtypea`\arrowtypeb`\arrowtypec;\height>}}

\def\dtrianglep<#1>[#2`#3`#4;#5`#6`#7]{{
\settriparms[#1]
\width=\height                         
\xext=\width                           
\yext=\height                          
\topadjust[#2``]
\botadjust[#3`#4`{#7}]
\leftsladjust[#3`#2`{#5}]
\rightadjust[#2`#4`{#6}]
\begin{picture}(\xext,\yext)(\xoff,\yoff)
\putdtrianglep<\arrowtypea`\arrowtypeb`\arrowtypec;\height>%
(0,0)[#2`#3`#4;#5`#6`{#7}]%
\end{picture}%
}}

\def\putbtrianglep<#1>(#2,#3)[#4`#5`#6;#7`#8`#9]{{%
\settriparms[#1]%
\xpos=#2 \ypos=#3
\puthmorphism(\xpos,\ypos)[#5`#6`{#9}]{\height}{\arrowtypec}b%
\advance\ypos by\height
\putmorphism(\xpos,\ypos)(1,-1)[``{#8}]{\height}{\arrowtypeb}r%
\putvmorphism(\xpos,\ypos)[#4``{#7}]{\height}{\arrowtypea}l%
}}

\def\putbtriangle{\@ifnextchar <{\putbtrianglep}{\putbtrianglep
   <\arrowtypea`\arrowtypeb`\arrowtypec;\height>}}
\def\btriangle{\@ifnextchar <{\btrianglep}{\btrianglep
   <\arrowtypea`\arrowtypeb`\arrowtypec;\height>}}

\def\btrianglep<#1>[#2`#3`#4;#5`#6`#7]{{
\settriparms[#1]
\width=\height                         
\xext=\width                           
\yext=\height                          
\topadjust[#2``]
\botadjust[#3`#4`{#7}]
\leftadjust[#2`#3`{#5}]
\rightsladjust[#4`#2`{#6}]
\begin{picture}(\xext,\yext)(\xoff,\yoff)
\putbtrianglep<\arrowtypea`\arrowtypeb`\arrowtypec;\height>%
(0,0)[#2`#3`#4;#5`#6`{#7}]%
\end{picture}%
}}

\def\putAtrianglep<#1>(#2,#3)[#4`#5`#6;#7`#8`#9]{{%
\settriparms[#1]%
\xpos=#2 \ypos=#3
{\multiply \height by2
\puthmorphism(\xpos,\ypos)[#5`#6`{#9}]{\height}{\arrowtypec}b}%
\advance\xpos by\height \advance\ypos by\height
\putmorphism(\xpos,\ypos)(-1,-1)[#4``{#7}]{\height}{\arrowtypea}l%
\putmorphism(\xpos,\ypos)(1,-1)[``{#8}]{\height}{\arrowtypeb}r%
}}

\def\putAtriangle{\@ifnextchar <{\putAtrianglep}{\putAtrianglep
   <\arrowtypea`\arrowtypeb`\arrowtypec;\height>}}
\def\Atriangle{\@ifnextchar <{\Atrianglep}{\Atrianglep
   <\arrowtypea`\arrowtypeb`\arrowtypec;\height>}}

\def\Atrianglep<#1>[#2`#3`#4;#5`#6`#7]{{
\settriparms[#1]
\width=\height                         
\xext=\width                           
\yext=\height                          
\topadjust[#2``]
\botadjust[#3`#4`{#7}]
\multiply \xext by2 
\leftsladjust[#3`#2`{#5}]
\rightsladjust[#4`#2`{#6}]
\begin{picture}(\xext,\yext)(\xoff,\yoff)%
\putAtrianglep<\arrowtypea`\arrowtypeb`\arrowtypec;\height>%
(0,0)[#2`#3`#4;#5`#6`{#7}]%
\end{picture}%
}}

\def\putAtrianglepairp<#1>(#2)[#3;#4`#5`#6`#7`#8]{{
\settripairparms[#1]%
\setpos(#2)%
\settokens[#3]%
\puthmorphism(\xpos,\ypos)[\tokenb`\tokenc`{#7}]{\height}{\arrowtyped}b%
\advance\xpos by\height
\advance\ypos by\height
\putmorphism(\xpos,\ypos)(-1,-1)[\tokena``{#4}]{\height}{\arrowtypea}l%
\putvmorphism(\xpos,\ypos)[``{#5}]{\height}{\arrowtypeb}m%
\putmorphism(\xpos,\ypos)(1,-1)[``{#6}]{\height}{\arrowtypec}r%
}}

\def\putAtrianglepair{\@ifnextchar <{\putAtrianglepairp}{\putAtrianglepairp%
   <\arrowtypea`\arrowtypeb`\arrowtypec`\arrowtyped`\arrowtypee;\height>}}
\def\Atrianglepair{\@ifnextchar <{\Atrianglepairp}{\Atrianglepairp%
   <\arrowtypea`\arrowtypeb`\arrowtypec`\arrowtyped`\arrowtypee;\height>}}

\def\Atrianglepairp<#1>[#2;#3`#4`#5`#6`#7]{{%
\settripairparms[#1]%
\settokens[#2]%
\width=\height
\xext=\width
\yext=\height
\topadjust[\tokena``]%
\vertadjust[\tokenb`\tokenc`{#6}]
\tempcountd=\tempcounta                       
\vertadjust[\tokenc`\tokend`{#7}]
\ifnum\tempcounta<\tempcountd                 
\tempcounta=\tempcountd\fi                    
\advance \yext by\tempcounta                  
\advance \yoff by-\tempcounta                 %
\multiply \xext by2 
\leftsladjust[\tokenb`\tokena`{#3}]
\rightsladjust[\tokend`\tokena`{#5}]%
\begin{picture}(\xext,\yext)(\xoff,\yoff)%
\putAtrianglepairp
<\arrowtypea`\arrowtypeb`\arrowtypec`\arrowtyped`\arrowtypee;\height>%
(0,0)[#2;#3`#4`#5`#6`{#7}]%
\end{picture}%
}}

\def\putVtrianglep<#1>(#2,#3)[#4`#5`#6;#7`#8`#9]{{%
\settriparms[#1]%
\xpos=#2 \ypos=#3
\advance\ypos by\height
{\multiply\height by2
\puthmorphism(\xpos,\ypos)[#4`#5`{#7}]{\height}{\arrowtypea}a}%
\putmorphism(\xpos,\ypos)(1,-1)[`#6`{#8}]{\height}{\arrowtypeb}l%
\advance\xpos by\height
\advance\xpos by\height
\putmorphism(\xpos,\ypos)(-1,-1)[``{#9}]{\height}{\arrowtypec}r%
}}

\def\putVtriangle{\@ifnextchar <{\putVtrianglep}{\putVtrianglep
   <\arrowtypea`\arrowtypeb`\arrowtypec;\height>}}
\def\Vtriangle{\@ifnextchar <{\Vtrianglep}{\Vtrianglep
   <\arrowtypea`\arrowtypeb`\arrowtypec;\height>}}

\def\Vtrianglep<#1>[#2`#3`#4;#5`#6`#7]{{
\settriparms[#1]
\width=\height                         
\xext=\width                           
\yext=\height                          
\topadjust[#2`#3`{#5}]
\botadjust[#4``]
\multiply \xext by2 
\leftsladjust[#2`#3`{#6}]
\rightsladjust[#3`#4`{#7}]
\begin{picture}(\xext,\yext)(\xoff,\yoff)%
\putVtrianglep<\arrowtypea`\arrowtypeb`\arrowtypec;\height>%
(0,0)[#2`#3`#4;#5`#6`{#7}]%
\end{picture}%
}}

\def\putVtrianglepairp<#1>(#2)[#3;#4`#5`#6`#7`#8]{{
\settripairparms[#1]%
\setpos(#2)%
\settokens[#3]%
\advance\ypos by\height
\putmorphism(\xpos,\ypos)(1,-1)[`\tokend`{#6}]{\height}{\arrowtypec}l%
\puthmorphism(\xpos,\ypos)[\tokena`\tokenb`{#4}]{\height}{\arrowtypea}a%
\advance\xpos by\height
\putvmorphism(\xpos,\ypos)[``{#7}]{\height}{\arrowtyped}m%
\advance\xpos by\height
\putmorphism(\xpos,\ypos)(-1,-1)[``{#8}]{\height}{\arrowtypee}r%
}}

\def\putVtrianglepair{\@ifnextchar <{\putVtrianglepairp}{\putVtrianglepairp%
    <\arrowtypea`\arrowtypeb`\arrowtypec`\arrowtyped`\arrowtypee;\height>}}
\def\Vtrianglepair{\@ifnextchar <{\Vtrianglepairp}{\Vtrianglepairp%
    <\arrowtypea`\arrowtypeb`\arrowtypec`\arrowtyped`\arrowtypee;\height>}}

\def\Vtrianglepairp<#1>[#2;#3`#4`#5`#6`#7]{{%
\settripairparms[#1]%
\settokens[#2]
\xext=\height                  
\width=\height                 
\yext=\height                  
\vertadjust[\tokena`\tokenb`{#4}]
\tempcountd=\tempcounta        
\vertadjust[\tokenb`\tokenc`{#5}]
\ifnum\tempcounta<\tempcountd%
\tempcounta=\tempcountd\fi
\advance \yext by\tempcounta
\botadjust[\tokend``]%
\multiply \xext by2
\leftsladjust[\tokena`\tokend`{#6}]%
\rightsladjust[\tokenc`\tokend`{#7}]%
\begin{picture}(\xext,\yext)(\xoff,\yoff)%
\putVtrianglepairp
<\arrowtypea`\arrowtypeb`\arrowtypec`\arrowtyped`\arrowtypee;\height>%
(0,0)[#2;#3`#4`#5`#6`{#7}]%
\end{picture}%
}}

\def\putCtrianglep<#1>(#2,#3)[#4`#5`#6;#7`#8`#9]{{%
\settriparms[#1]%
\xpos=#2 \ypos=#3
\advance\ypos by\height
\putmorphism(\xpos,\ypos)(1,-1)[``{#9}]{\height}{\arrowtypec}l%
\advance\xpos by\height
\advance\ypos by\height
\putmorphism(\xpos,\ypos)(-1,-1)[#4`#5`{#7}]{\height}{\arrowtypea}l%
{\multiply\height by 2
\putvmorphism(\xpos,\ypos)[`#6`{#8}]{\height}{\arrowtypeb}r}%
}}

\def\putCtriangle{\@ifnextchar <{\putCtrianglep}{\putCtrianglep
    <\arrowtypea`\arrowtypeb`\arrowtypec;\height>}}
\def\Ctriangle{\@ifnextchar <{\Ctrianglep}{\Ctrianglep
    <\arrowtypea`\arrowtypeb`\arrowtypec;\height>}}

\def\Ctrianglep<#1>[#2`#3`#4;#5`#6`#7]{{
\settriparms[#1]
\width=\height                          
\xext=\width                            
\yext=\height                           
\multiply \yext by2 
\topadjust[#2``]
\botadjust[#4``]
\sladjust[#3`#2`{#5}]{\width}
\tempcountd=\tempcounta                 
\sladjust[#3`#4`{#7}]{\width}
\ifnum \tempcounta<\tempcountd          
\tempcounta=\tempcountd\fi              
\advance \xext by\tempcounta            
\advance \xoff by-\tempcounta           %
\rightadjust[#2`#4`{#6}]
\begin{picture}(\xext,\yext)(\xoff,\yoff)%
\putCtrianglep<\arrowtypea`\arrowtypeb`\arrowtypec;\height>%
(0,0)[#2`#3`#4;#5`#6`{#7}]%
\end{picture}%
}}

\def\putDtrianglep<#1>(#2,#3)[#4`#5`#6;#7`#8`#9]{{%
\settriparms[#1]%
\xpos=#2 \ypos=#3
\advance\xpos by\height \advance\ypos by\height
\putmorphism(\xpos,\ypos)(-1,-1)[``{#9}]{\height}{\arrowtypec}r%
\advance\xpos by-\height \advance\ypos by\height
\putmorphism(\xpos,\ypos)(1,-1)[`#5`{#8}]{\height}{\arrowtypeb}r%
{\multiply\height by 2
\putvmorphism(\xpos,\ypos)[#4`#6`{#7}]{\height}{\arrowtypea}l}%
}}

\def\putDtriangle{\@ifnextchar <{\putDtrianglep}{\putDtrianglep
    <\arrowtypea`\arrowtypeb`\arrowtypec;\height>}}
\def\Dtriangle{\@ifnextchar <{\Dtrianglep}{\Dtrianglep
   <\arrowtypea`\arrowtypeb`\arrowtypec;\height>}}

\def\Dtrianglep<#1>[#2`#3`#4;#5`#6`#7]{{
\settriparms[#1]
\width=\height                         
\xext=\height                          
\yext=\height                          
\multiply \yext by2 
\topadjust[#2``]
\botadjust[#4``]
\leftadjust[#2`#4`{#5}]
\sladjust[#3`#2`{#5}]{\height}
\tempcountd=\tempcountd                
\sladjust[#3`#4`{#7}]{\height}
\ifnum \tempcounta<\tempcountd         
\tempcounta=\tempcountd\fi             
\advance \xext by\tempcounta           %
\begin{picture}(\xext,\yext)(\xoff,\yoff)
\putDtrianglep<\arrowtypea`\arrowtypeb`\arrowtypec;\height>%
(0,0)[#2`#3`#4;#5`#6`{#7}]%
\end{picture}%
}}

\def\setrecparms[#1`#2]{\width=#1 \height=#2}%
\def\recursep<#1`#2>[#3;#4`#5`#6`#7`#8]{{%
\width=#1 \height=#2
\settokens[#3]
\settowidth{\tempdimen}{$\tokena$}
\ifdim\tempdimen=0pt
  \savebox{\tempboxa}{\hbox{$\tokenb$}}%
  \savebox{\tempboxb}{\hbox{$\tokend$}}%
  \savebox{\tempboxc}{\hbox{$#6$}}%
\else
  \savebox{\tempboxa}{\hbox{$\hbox{$\tokena$}\times\hbox{$\tokenb$}$}}%
  \savebox{\tempboxb}{\hbox{$\hbox{$\tokena$}\times\hbox{$\tokend$}$}}%
  \savebox{\tempboxc}{\hbox{$\hbox{$\tokena$}\times\hbox{$#6$}$}}%
\fi
\ypos=\height
\divide\ypos by 2
\xpos=\ypos
\advance\xpos by \width
\xext=\xpos \yext=\height
\topadjust[#3`\usebox{\tempboxa}`{#4}]%
\botadjust[#5`\usebox{\tempboxb}`{#8}]%
\sladjust[\tokenc`\tokenb`{#5}]{\ypos}%
\tempcountd=\tempcounta
\sladjust[\tokenc`\tokend`{#5}]{\ypos}%
\ifnum \tempcounta<\tempcountd
\tempcounta=\tempcountd\fi
\advance \xext by\tempcounta
\advance \xoff by-\tempcounta
\rightadjust[\usebox{\tempboxa}`\usebox{\tempboxb}`\usebox{\tempboxc}]%
\bfig
\putCtrianglep<-1`1`1;\ypos>(0,0)[`\tokenc`;#5`#6`{#7}]%
\puthmorphism(\ypos,0)[\tokend`\usebox{\tempboxb}`{#8}]{\width}{-1}b%
\puthmorphism(\ypos,\height)[\tokenb`\usebox{\tempboxa}`{#4}]{\width}{-1}a%
\advance\ypos by \width
\putvmorphism(\ypos,\height)[``\usebox{\tempboxc}]{\height}1r%
\efig
}}

\def\recurse{\@ifnextchar <{\recursep}{\recursep<\width`\height>}}

\def\puttwohmorphisms(#1,#2)[#3`#4;#5`#6]#7#8#9{{%
\puthmorphism(#1,#2)[#3`#4`]{#7}0a
\ypos=#2
\advance\ypos by 20
\puthmorphism(#1,\ypos)[\phantom{#3}`\phantom{#4}`#5]{#7}{#8}a
\advance\ypos by -40
\puthmorphism(#1,\ypos)[\phantom{#3}`\phantom{#4}`#6]{#7}{#9}b
}}

\def\puttwovmorphisms(#1,#2)[#3`#4;#5`#6]#7#8#9{{%
\putvmorphism(#1,#2)[#3`#4`]{#7}0a
\xpos=#1
\advance\xpos by -20
\putvmorphism(\xpos,#2)[\phantom{#3}`\phantom{#4}`#5]{#7}{#8}l
\advance\xpos by 40
\putvmorphism(\xpos,#2)[\phantom{#3}`\phantom{#4}`#6]{#7}{#9}r
}}

\def\puthcoequalizer(#1)[#2`#3`#4;#5`#6`#7]#8#9{{%
\setpos(#1)%
\puttwohmorphisms(\xpos,\ypos)[#2`#3;#5`#6]{#8}11%
\advance\xpos by #8
\puthmorphism(\xpos,\ypos)[\phantom{#3}`#4`#7]{#8}1{#9}
}}

\def\putvcoequalizer(#1)[#2`#3`#4;#5`#6`#7]#8#9{{%
\setpos(#1)%
\puttwovmorphisms(\xpos,\ypos)[#2`#3;#5`#6]{#8}11%
\advance\ypos by -#8
\putvmorphism(\xpos,\ypos)[\phantom{#3}`#4`#7]{#8}1{#9}
}}

\def\putthreehmorphisms(#1)[#2`#3;#4`#5`#6]#7(#8)#9{{%
\setpos(#1) \settypes(#8)
\if a#9 %
     \vertsize{\tempcounta}{#5}%
     \vertsize{\tempcountb}{#6}%
     \ifnum \tempcounta<\tempcountb \tempcounta=\tempcountb \fi
\else
     \vertsize{\tempcounta}{#4}%
     \vertsize{\tempcountb}{#5}%
     \ifnum \tempcounta<\tempcountb \tempcounta=\tempcountb \fi
\fi
\advance \tempcounta by 60
\puthmorphism(\xpos,\ypos)[#2`#3`#5]{#7}{\arrowtypeb}{#9}
\advance\ypos by \tempcounta
\puthmorphism(\xpos,\ypos)[\phantom{#2}`\phantom{#3}`#4]{#7}{\arrowtypea}{#9}
\advance\ypos by -\tempcounta \advance\ypos by -\tempcounta
\puthmorphism(\xpos,\ypos)[\phantom{#2}`\phantom{#3}`#6]{#7}{\arrowtypec}{#9}
}}

\def\putarc(#1,#2)[#3`#4`#5]#6#7#8{{%
\xpos #1
\ypos #2
\width #6
\arrowlength #6
\putbox(\xpos,\ypos){#3\vphantom{#4}}%
{\advance \xpos by\arrowlength
\putbox(\xpos,\ypos){\vphantom{#3}#4}}%
\horsize{\tempcounta}{#3}%
\horsize{\tempcountb}{#4}%
\divide \tempcounta by2
\divide \tempcountb by2
\advance \tempcounta by30
\advance \tempcountb by30
\advance \xpos by\tempcounta
\advance \arrowlength by-\tempcounta
\advance \arrowlength by-\tempcountb
\halflength=\arrowlength \divide\halflength by 2
\divide\arrowlength by 5
\put(\xpos,\ypos){\bezier{\arrowlength}(0,0)(50,50)(\halflength,50)}
\ifnum #7=-1 \put(\xpos,\ypos){\vector(-3,-2)0} \fi
\advance\xpos by \halflength
\put(\xpos,\ypos){\xpos=\halflength \advance\xpos by -50
   \bezier{\arrowlength}(0,50)(\xpos,50)(\halflength,0)}
\ifnum #7=1 {\advance \xpos by
   \halflength \put(\xpos,\ypos){\vector(3,-2)0}} \fi
\advance\ypos by 50
\vertsize{\tempcounta}{#5}%
\divide\tempcounta by2
\advance \tempcounta by20
\if a#8 %
   \advance \ypos by\tempcounta
   \putbox(\xpos,\ypos){#5}%
\else
   \advance \ypos by-\tempcounta
   \putbox(\xpos,\ypos){#5}%
\fi
}}

\makeatother

\usepackage{amsthm}
\usepackage{amsmath}
\usepackage{dsfont}
\usepackage{stmaryrd}
\usepackage{graphics}
\usepackage[linktocpage]{hyperref}

\newtheorem{theorem}{Theorem}[section]
\newtheorem{lemma}[theorem]{Lemma}
\newtheorem{corollary}[theorem]{Corollary}

\newtheorem{proposition}[theorem]{Proposition}

\makeindex \makeglossary
\begin{document}

\sloppy

\newcommand{\nl}{\hspace{2cm}\\ }

\def\nec{\Box}
\def\pos{\Diamond}
\def\diam{{\tiny\Diamond}}

\def\lc{\lceil}
\def\rc{\rceil}
\def\lf{\lfloor}
\def\rf{\rfloor}
\def\lk{\langle}
\def\rk{\rangle}
\def\blk{\dot{\langle\!\!\langle}}
\def\brk{\dot{\rangle\!\!\rangle}}

\newcommand{\pa}{\parallel}
\newcommand{\lra}{\longrightarrow}
\newcommand{\hra}{\hookrightarrow}
\newcommand{\hla}{\hookleftarrow}
\newcommand{\ra}{\rightarrow}
\newcommand{\la}{\leftarrow}
\newcommand{\lla}{\longleftarrow}
\newcommand{\da}{\downarrow}
\newcommand{\ua}{\uparrow}
\newcommand{\dA}{\downarrow\!\!\!^\bullet}
\newcommand{\uA}{\uparrow\!\!\!_\bullet}
\newcommand{\Da}{\Downarrow}
\newcommand{\DA}{\Downarrow\!\!\!^\bullet}
\newcommand{\UA}{\Uparrow\!\!\!_\bullet}
\newcommand{\Ua}{\Uparrow}
\newcommand{\Lra}{\Longrightarrow}
\newcommand{\Ra}{\Rightarrow}
\newcommand{\Lla}{\Longleftarrow}
\newcommand{\La}{\Leftarrow}
\newcommand{\nperp}{\perp\!\!\!\!\!\setminus\;\;}
\newcommand{\pq}{\preceq}

\newcommand{\lms}{\longmapsto}
\newcommand{\ms}{\mapsto}
\newcommand{\subseteqnot}{\subseteq\hskip-4 mm_\not\hskip3 mm}

\def\o{{\omega}}
\def\sM{{\bf sM}}
\def\bA{{\bf A}}
\def\bEM{{\bf EM}}
\def\bM{{\bf M}}
\def\bN{{\bf N}}
\def\bC{{\bf C}}
\def\bI{{\bf I}}
\def\bK{{\bf K}}
\def\bL{{\bf L}}
\def\bT{{\bf T}}
\def\bS{{\bf S}}
\def\bD{{\bf D}}
\def\bB{{\bf B}}
\def\bW{{\bf W}}
\def\bP{{\bf P}}
\def\bX{{\bf X}}
\def\bU{{\bf U}}
\def\bY{{\bf Y}}
\def\ba{{\bf a}}
\def\bb{{\bf b}}
\def\bc{{\bf c}}
\def\bd{{\bf d}}
\def\bh{{\bf h}}
\def\bi{{\bf i}}
\def\bj{{\bf j}}
\def\bk{{\bf k}}
\def\bm{{\bf m}}
\def\bn{{\bf n}}
\def\bp{{\bf p}}
\def\bq{{\bf q}}
\def\be{{\bf e}}
\def\br{{\bf r}}
\def\bi{{\bf i}}
\def\bs{{\bf s}}
\def\bt{{\bf t}}
\def\bu{{\bf u}}
\def\bv{{\bf v}}
\def\bw{{\bf w}}
\def\bz{{\bf z}}

\def\zero{{\bf 0}}
\def\jeden{{\bf 1}}
\def\dwa{{\bf 2}}
\def\trzy{{\bf 3}}
\def\Lam{{\bf \Lambda}}
\def\Tree{{\bf Tree}}
\def\cTree{{\bf Zoom}}
\def\tZoom{{\bf Zoom_t}}

\def\CS{{\bf CS}}
\def\pZoom{{\bf pZoom}}
\def\cTree{{\bf pZoom}}
\def\wZoom{{\bf wZoom}}

\def\cBL{{\cal BL}}
\def\cB{{\cal B}}
\def\cA{{\cal A}}
\def\cC{{\cal C}}
\def\cD{{\cal D}}
\def\cE{{\cal E}}
\def\cEM{{\cal EM}}
\def\cF{{\cal F}}
\def\cG{{\cal G}}
\def\cH{{\cal H}}
\def\cI{{\cal I}}
\def\cJ{{\cal J}}
\def\cK{{\cal K}}
\def\cL{{\cal L}}
\def\cN{{\cal N}}
\def\cM{{\cal M}}
\def\cO{{\cal O}}
\def\cP{{\cal P}}
\def\cQ{{\cal Q}}
\def\cR{{\cal R}}
\def\cS{{\cal S}}
\def\cT{{\cal T}}
\def\cU{{\cal U}}
\def\cV{{\cal V}}
\def\cW{{\cal W}}
\def\cX{{\cal X}}
\def\cY{{\cal Y}}

\def\sgn{{\bf sgn}}
\def\el{{\bf el}}
\def\ll{{\bf ll}}
\def\nextF{{\bf nextF}}
\def\prevF{{\bf prevF}}
\def\nextP{{\bf nextP}}
\def\prevP{{\bf prevP}}
{\def\next{{\bf next}}}
\def\neig{{\bf neig}}
\def\neigl{{\bf neig}_{low}}
\def\neigh{{\bf neig}_{high}}
\def\next{{\bf succ}}
\def\prev{{\bf prev}}
\def\top{{\bf top}}
\def\tr{{\bf tr}}
\def\pu{{\bf pu}}
\def\Fl{{\bf Fl}}
\def\Pl{{\bf Pl}}
\def\pFl{{\bf pFl}}
\def\Flags{{\bf Flags}}
\def\Path{{\bf Path}}
\def\diff{{\bf diff}}
\def\id{{\bf id}}
\def\pFlags{{\bf pFlags}}
\def\MpFlags{{\bf MpFlags}}
\def\MFlags{{\bf MFlags}}
\def\Cyl{{\bf Cyl}}
\def\Cyli{{\bf Cyl}_\iota}
\def\Cylh{{\bf Cyl}_h}
\def\Cylhi{{\bf Cyl}_{h,\iota}}
\def\Cylp{{\bf Cyl}_p}
\def\Ope{{\bf Ope}}
\def\pOpe{{\bf pOpe}}
\def\pOpeCard{{\bf pOpeCard}}
\def\Opei{{\bf Ope}_\iota}
\def\pOpei{{\bf pOpe}_\iota}
\def\pOpeie{{\bf pOpe}_{\iota e}}
\def\Opeo{{\bf Ope}_\omega}
\def\pOpeo{{\bf pOpe}_\omega}
\def\pHg{{\bf pHg}}
\def\pHgi{{\bf pHg_{\iota}}}
\def\pOHgi{{\bf pOHg_{\iota}}}
\def\cHi{{\cal H}_{\iota}}
\def\Fp{{\bf Fp}}
\def\sgn{{\bf sgn}}
\def\lhds{{\lhd\!\!\!\cdot\;\,}}
\def\inn{{\in\!\!\!\rightharpoondown}}
\def\St{{\bf St}}
\def\pOpe{{\bf pOpe}}
\def\pOpem{{\bf pOpem}}
\def\pOpeCard{{\bf pOpeCard}}
\def\cThick{{\bf tZoom}}
\def\tZoom{{\bf tZoom}}
\def\OO{{\bf OpeOpe}}
\def\cvr{ \rm{cvr}}
\def\lvs{ {\rm{lvs}}}
\def\sup{ \rm{sup}}

\def\sdiff{\stackrel{.}{-}}

\def\LMF{{\bf LMF}}
\def\Mon{{\bf Mon}}
\def\LAdj{{\bf LAdj}}
\def\Adj{{\bf Adj}}
\def\RAdj{{\bf RAdj}}
\def\Act{{\bf Act}}
\def\Clsd{{\bf Closed}}
\def\Mlcv{{\bf Malcev}}
\def\cMlcv{{\bf coMalcev}}
\def\Mnd{{\bf Mnd}}
\def\bCat{{{\bf Cat}}}
\def\Cat{{{\bf Cat}}}
\def\Gpd{{{\bf Gpd}}}
\def\Gph{{{\bf Gph}}}
\def\Fib{{{\bf Fib}}}
\def\DFib{{{\bf DFib}}}
\def\BiFib{{{\bf BiFib}}}
\def\Catrc{{{\bf Cat}_{rc}}}
\def\Monrc{{{\bf Mon}_{rc}}}
\def\Mod{{\bf Mod}}
\def\Modm{{\bf cMod}}
\def\cMod{{\bf cMod}}
\def\mon{{{\bf mon}}}
\def\act{{{\bf act}}}
\def\bem{{{\bf em}}}
\def\bkem{{{\bf kem}}}
\def\bKl{{{\bf Kl}}}

\def\oG{{{\omega}Gr}}
\def\mts{{MltSet}}
\def\pOpeSet{{\bf pOpeSet}}


\def\oCat{{\bf \o Cat}}
\def\oC{{\bf \o C}}
\def\oGph{{\bf \o Gph}}
\def\AMon{{\bf AnMnd}}
\def\An{{\bf An}}
\def\Poly{{\bf Poly}}
\def\San{{\bf San}}
\def\Taut{{\bf Taut}}
\def\PMnd{{\bf PolyMnd}}
\def\SanMnd{{\bf SanMnd}}
\def\RiMnd{{\bf RiMnd}}
\def\End{{\bf End}}

\def\Diag{{\bf Diag}}

\def\ET{\bf ET}
\def\RegET{\bf RegET}
\def\RET{\bf RegET}
\def\LrET{\bf LrET}
\def\RiET{\bf RiET}
\def\SregET{\bf SregET}
\def\Cart{\bf Cart}
\def\wCart{\bf wCart}
\def\CartMnd{\bf CartMnd}
\def\wCartMnd{\bf wCartMnd}

\def\LT{\bf LT}
\def\RegLT{\bf RegLT}
\def\ALT{\bf AnLT}
\def\RiLT{\bf RiLT}

\def\FOp{\bf FOp}
\def\RegOp{\bf RegOp}
\def\SOp{\bf SOp}
\def\RiOp{\bf RiOp}

\def\MonCat{{{\bf MonCat}}}
\def\ActMonCat{{{\bf ActMonCat}}}

\def\F{\mathds{F}}
\def\E{\mathds{E}}
\def\S{\mathds{S}}
\def\I{\mathds{I}}
\def\B{\mathds{B}}

\def\V{\mathds{V}}
\def\W{\mathds{W}}
\def\M{\mathds{M}}
\def\N{\mathds{N}}
\def\R{\mathds{R}}

\pagenumbering{arabic} \setcounter{page}{1}

\title{\bf\Large Duality for positive opetopes and positive zoom complexes}

\author{Marek Zawadowski\\
}
\maketitle

\begin{abstract} We show that the (positive) zoom complexes, with fairly natural morphisms, form a dual category to the category of positive opetopes with contraction epimorphisms. We also show how this duality can be slightly extended to opetopic cardinals.
\end{abstract}

\tableofcontents

\section{Introduction}

The opetopes are higher dimensional shapes that were originally invented in [BD] as shapes that can be used to define a notion of higher dimensional category. By now there are more than a dozen of definitions of opetopes. Some definitions use very abstract categorical machinery \cite{Bu}, \cite{Z3}, some are more concrete using one way or another some kinds of operads and/or polynomial or analytic monads  \cite{BD}, \cite{HMP}, \cite{Lei}, \cite{Cheng}, \cite{Z3}, \cite{SZ}, \cite{FS}, some definitions describe the ways opetopes can be generated \cite{HMP}, \cite{CTM}, \cite{T},  and finally there are also some purely combinatorial definitions \cite{Palm}, \cite{KJBM}, \cite{Z1}, \cite{Z2}, \cite{St}.

So it is not surprising that it is easier to show a picture of an opetope than to give a simple definition that will leave the reader with no doubts as to what an opetope is.
In this paper we will deal with positive opetopes only. This means that each face in opetopes we consider has at least one face of codimension $1$ in its domain. As a consequence such opetopes have no loops.

The definitions of opetopes mentioned above seem to agree in a `reasonable sense'. However, the morphisms between opetopes are not always treated the same way. In some approaches even face maps between opetopes do not seem to be natural. In \cite{Z1},\cite{Z3} it is shown how opetopes can be treated as some special kinds of $\o$-categories and therefore all the $\o$-functors (i.e. all face maps and all degeneracies) between opetopes can be considered. Of all the definitions of opetopes,  the one given in \cite{KJBM}, through so called zoom complexes, seem to be very different from any other. In fact, to describe this definition even pictures are of a very different kind, as the reader can notice below. In this paper we shall show that there is an explanation of this phenomenon. Namely, the (positive) zoom complexes, with fairly natural morphisms, form a dual category to the category of positive opetopes with contraction epimorphisms (often called $\iota$-epis).  The contraction epimorphisms are some kind of degeneracies of opetopes that can send a face only to a face but possibly of a lower dimension (still preserving usual constraints concerning both domains and codomains).

Below we draw some pictures of positive opetopes of few low dimensions and corresponding dual (positive) zoom complexes, a simplified version of zoom complexes introduced in \cite{KJBM}.

An opetope $O_1$ of dimension 1
\begin{center} \xext=500 \yext=150
\begin{picture}(\xext,\yext)(\xoff,\yoff)
 \putmorphism(0,50)(1,0)[t_2`t_1`y]{500}{1}a
 \end{picture}
\end{center}
and its dual positive zoom complex $T_1$
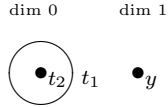
\begin{center} \xext=500 \yext=400
\begin{picture}(\xext,\yext)(\xoff,\yoff)

 \put(40,350){\tiny dim 0}
  \put(420,350){\tiny dim 1}

 \put(150,150){\circle{220}}
 \put(125,125){$\bullet _{t_2}$}
 \put(290,125){$ _{t_1}$}
 \put(460,125){$\bullet _y$}
 \end{picture}
\end{center}

An opetope $O_2$ of dimension 2
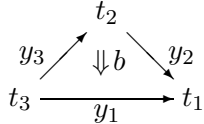
\begin{center} \xext=500 \yext=400
\begin{picture}(\xext,\yext)(\xoff,\yoff)
 \settriparms[-1`1`1;300]
 \putAtriangle(50,50)[t_2`t_3`t_1;y_3`y_2`y_1]
 \put(300,150){$\Da\! b$}
 \end{picture}
\end{center}
and its dual positive zoom complex $T_2$
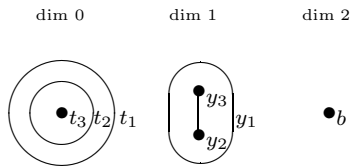
\begin{center} \xext=750 \yext=550
\begin{picture}(\xext,\yext)(\xoff,\yoff)
 \put(60,470){\tiny dim 0}
  \put(520,470){\tiny dim 1}
  \put(980,470){\tiny dim 2}
 \put(150,150){\circle{220}}
 \put(150,150){\circle{350}}
 \put(125,125){$\bullet _{t_3}$}
 \put(260,125){$ _{t_2}$}
 \put(350,125){$ _{t_1}$}

 \put(600,50){$\bullet _{y_2}$}
 \put(600,200){$\bullet _{y_3}$}
 \put(620,70){\line(0,1){150}}
 \put(630,150){\oval(220,350)}
 \put(750,125){$_{y_1}$}

 \put(1050,125){$\bullet _b$}

 \end{picture}
\end{center}

An opetope $O_3$ of dimension 3
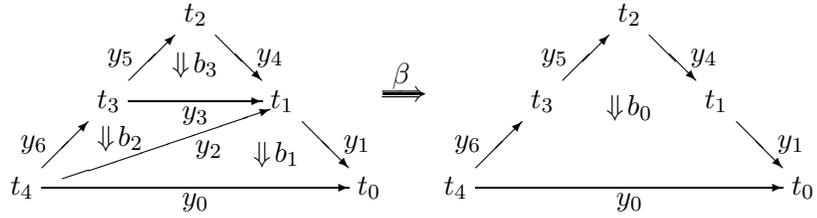
\begin{center} \xext=2400 \yext=800
\begin{picture}(\xext,\yext)(\xoff,\yoff)
 \label{figure ppof}
 \settriparms[-1`1`1;300]
 \putAtriangle(300,350)[t_2`t_3`t_1;y_5`y_4`y_3]
 \put(520,450){$\Da\! b_3$}
 \settriparms[-1`0`0;300]
 \putAtriangle(0,50)[\phantom{t_3}`t_4`;y_6``]
 \settriparms[0`1`0;300]
 \putAtriangle(600,50)[\phantom{t_1}``t_0;`y_1`]
 \putmorphism(0,50)(1,0)[\phantom{t_4}`\phantom{t_0}`y_0]{1200}{1}b
 \putmorphism(350,150)(3,1)[\phantom{t_4}`\phantom{t_0}`]{300}{1}a
 \put(600,170){$y_2$}
 \put(260,195){$\Da\! b_2$}
 \put(800,140){$\Da\! b_1$}

\put(1240,350){$\Longrightarrow$}
  \put(1250,375){\line(1,0){135}}
   \put(1280,410){$\beta$}

\settriparms[-1`1`0;300]
 \putAtriangle(1800,350)[t_2`t_3`t_1;y_5`y_4`]
 \put(2020,300){$\Da\! b_0$}
 \settriparms[-1`0`0;300]
 \putAtriangle(1500,50)[\phantom{t_3}`t_4`;y_6``]
 \settriparms[0`1`0;300]
 \putAtriangle(2100,50)[\phantom{t_1}``t_0;`y_1`]
 \putmorphism(1500,50)(1,0)[\phantom{t_4}`\phantom{t_0}`y_0]{1200}{1}b
\end{picture}
\end{center}
and its dual positive zoom complex $T_3$
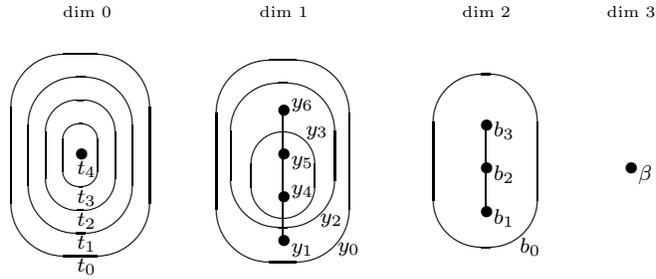
\begin{center} \xext=2200 \yext=950
\begin{picture}(\xext,\yext)(\xoff,\yoff)
 \put(260,900){\tiny dim 0}
  \put(940,900){\tiny dim 1}
  \put(1640,900){\tiny dim 2}
  \put(2140,900){\tiny dim 3}
 \put(300,400){$\bullet$}
 \put(310,370){$ _{t_4}$}
  \put(320,420){\oval(120,220)}
  \put(310,270){$ _{t_3}$}
  \put(320,420){\oval(240,380)}
  \put(310,190){$ _{t_2}$}
   \put(320,420){\oval(360,540)}
  \put(310,110){$ _{t_1}$}
   \put(320,420){\oval(480,700)}
  \put(310,30){$ _{t_0}$}

  \put(1000,100){$\bullet _{y_1}$}
 \put(1000,250){$\bullet ^{y_4}$}
  \put(1000,400){$\bullet _{y_5}$}
 \put(1000,550){$\bullet ^{y_6}$}
 \put(1020,120){\line(0,1){450}}

 \put(1020,350){\oval(220,300)}
 \put(1100,505){$ _{y_3}$}

 \put(1020,420){\oval(360,500)}
 \put(1150,200){$ _{y_2}$}

  \put(1020,400){\oval(460,700)}
 \put(1210,100){$ _{y_0}$}

  \put(1700,200){$\bullet _{b_1}$}
 \put(1700,350){$\bullet _{b_2}$}
  \put(1700,500){$\bullet _{b_3}$}
 \put(1720,220){\line(0,1){300}}

   \put(1720,400){\oval(360,600)}
 \put(1840,100){$ _{b_0}$}

 \put(2200,350){$\bullet _{\beta}$}

 \end{picture}
\end{center}
The drawings of the above positive zoom complexes are, in fact, drawings of consecutive (non-empty) constellations of positive zoom complexes. In particular, it is not an accident that the partial order of nesting of circles in one constellation is isomorphic to the partial order vertices of the next constellation.

Opetopic cardinals still consist of cells that can be meaningfully composed in a unique way, but have a bit more general shape. Here is an example
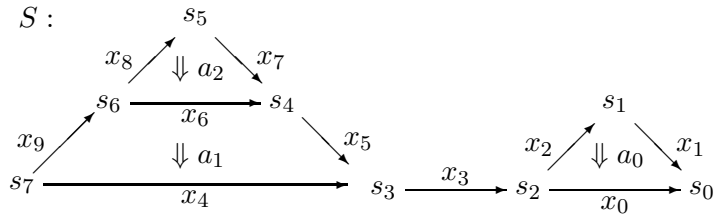
\begin{center} \xext=2400 \yext=750
\begin{picture}(\xext,\yext)(\xoff,\yoff)
 \settriparms[-1`1`1;300]
 \putAtriangle(300,350)[s_5`s_6`s_4;x_8`x_7`x_6]
 \put(560,420){\makebox(100,100){$\Downarrow a_2$}}
 \put(560,120){\makebox(100,100){$\Downarrow a_1$}}
 \putmorphism(0,70)(1,0)[s_7`\phantom{s_7}`x_4]{1200}{1}b
 \putmorphism(110,180)(1,1)[\phantom{s_7}`\phantom{s_3}`x_9]{70}{1}l
 \putmorphism(900,370)(1,-1)[\phantom{s_7}`\phantom{s_3}`x_5]{300}{1}r
 \putAtriangle(1750,50)[s_1`s_2`s_0;x_2`x_1`x_0]
 \put(2010,120){\makebox(100,100){$\Downarrow a_0$}}
 \put(0,600){\makebox(100,100){$S:$}}
\putmorphism(1250,50)(1,0)[s_3`\phantom{s_2}`x_3]{500}{1}a
\end{picture}
\end{center}
and its dual wide positive zoom complex

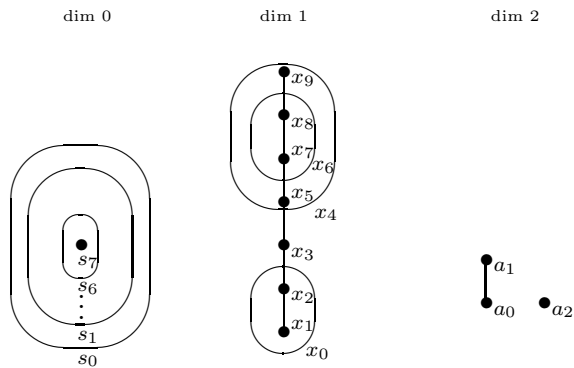
\begin{center} \xext=2200 \yext=1250
\begin{picture}(\xext,\yext)(\xoff,\yoff)
 \put(260,1200){\tiny dim 0}
  \put(940,1200){\tiny dim 1}
  \put(1740,1200){\tiny dim 2}

 \put(300,400){$\bullet$}
 \put(310,370){$ _{s_7}$}
  \put(320,420){\oval(120,220)}
  \put(310,280){$ _{s_6}$}
  \put(310,170){$\vdots$}
   \put(320,420){\oval(360,540)}
  \put(310,110){$ _{s_1}$}
   \put(320,420){\oval(480,700)}
  \put(310,30){$ _{s_0}$}

  \put(1000,100){$\bullet ^{x_1}$}
 \put(1000,250){$\bullet _{x_2}$}
  \put(1000,400){$\bullet _{x_3}$}
 \put(1000,550){$\bullet ^{x_5}$}
   \put(1000,700){$\bullet ^{x_7}$}
 \put(1000,850){$\bullet _{x_8}$}
   \put(1000,1000){$\bullet _{x_9}$}

 \put(1025,120){\line(0,1){900}}

 \put(1020,200){\oval(220,300)}
 \put(1100,55){$ _{x_0}$}

 \put(1020,800){\oval(220,300)}
 \put(1130,535){$ _{x_4}$}

 \put(1020,800){\oval(360,500)}
 \put(1120,700){$ _{x_6}$}

  \put(1700,200){$\bullet _{a_1}$}
 \put(1700,350){$\bullet _{a_2}$}
  \put(1900,200){$\bullet _{a_0}$}
 \put(1720,220){\line(0,1){150}}
 \end{picture}
\end{center}
To illustrate how the duality works on morphisms is a bit more involved. We shall present a $\iota$-epimorphism from $O_3$ to $O_2$ by naming faces of $O_3$
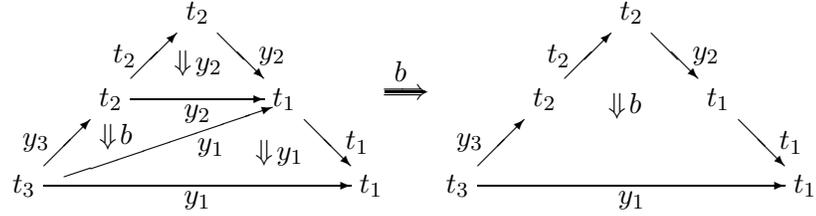
\begin{center} \xext=2400 \yext=800
\begin{picture}(\xext,\yext)(\xoff,\yoff)

 \settriparms[-1`1`1;300]
 \putAtriangle(300,350)[t_2`t_2`t_1;t_2`y_2`y_2]
 \put(520,450){$\Da\! y_2$}
 \settriparms[-1`0`0;300]
 \putAtriangle(0,50)[\phantom{t_3}`t_3`;y_3``]
 \settriparms[0`1`0;300]
 \putAtriangle(600,50)[\phantom{t_1}``t_1;`t_1`]
 \putmorphism(0,50)(1,0)[\phantom{t_4}`\phantom{t_0}`y_1]{1200}{1}b
 \putmorphism(350,150)(3,1)[\phantom{t_4}`\phantom{t_0}`]{300}{1}a
 \put(600,170){$y_1$}
 \put(260,195){$\Da\! b$}
 \put(800,140){$\Da\! y_1$}

\put(1240,350){$\Longrightarrow$}
  \put(1250,375){\line(1,0){135}}
   \put(1280,410){$b$}

\settriparms[-1`1`0;300]
 \putAtriangle(1800,350)[t_2`t_2`t_1;t_2`y_2`]
 \put(2020,300){$\Da\! b$}
 \settriparms[-1`0`0;300]
 \putAtriangle(1500,50)[\phantom{t_3}`t_3`;y_3``]
 \settriparms[0`1`0;300]
 \putAtriangle(2100,50)[\phantom{t_1}``t_1;`t_1`]
 \putmorphism(1500,50)(1,0)[\phantom{t_4}`\phantom{t_0}`y_1]{1200}{1}b
\end{picture}
\end{center}
by the cells on $O_2$
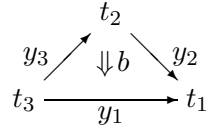
\begin{center} \xext=500 \yext=400
\begin{picture}(\xext,\yext)(\xoff,\yoff)

 \settriparms[-1`1`1;300]
 \putAtriangle(50,50)[t_2`t_3`t_1;y_3`y_2`y_1]
 \put(300,150){$\Da\! b$}

 \end{picture}
\end{center}
they are sent to. For example, there are three faces in $O_3$ sent to the $1$-face $y_2$. Two of them are $1$-faces and one of them is a $2$-face. The dual of this morphism is a collection of three embeddings of trees sending leaves (vertices) to leaves, and inner nodes (circles) to inner nodes and respecting the relation between circles at one dimension and the vertices on the next dimension. We present the dual morphism of positive zoom complexes from $T_2$ to $T_3$ drawing $T_3$ and naming its nodes
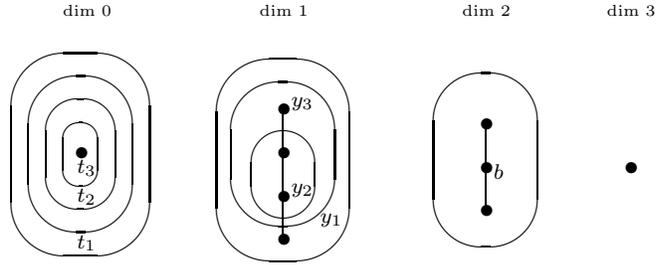
\begin{center} \xext=2200 \yext=950
\begin{picture}(\xext,\yext)(\xoff,\yoff)
 \put(260,900){\tiny dim 0}
  \put(940,900){\tiny dim 1}
  \put(1640,900){\tiny dim 2}
  \put(2140,900){\tiny dim 3}
 \put(300,400){$\bullet$}
 \put(310,370){$ _{t_3}$}
  \put(320,420){\oval(120,220)}
  \put(310,270){$ _{t_2}$}
  \put(320,420){\oval(240,380)}
   \put(320,420){\oval(360,540)}
  \put(310,110){$ _{t_1}$}
   \put(320,420){\oval(480,700)}

  \put(1000,100){$\bullet$}
 \put(1000,250){$\bullet ^{y_2}$}
  \put(1000,400){$\bullet $}
 \put(1000,550){$\bullet ^{y_3}$}
 \put(1020,120){\line(0,1){450}}

 \put(1020,350){\oval(220,300)}

 \put(1020,420){\oval(360,500)}
 \put(1150,200){$ _{y_1}$}

  \put(1020,400){\oval(460,700)}

  \put(1700,200){$\bullet$}
 \put(1700,350){$\bullet _{b}$}
  \put(1700,500){$\bullet$}
 \put(1720,220){\line(0,1){300}}

   \put(1720,400){\oval(360,600)}

 \put(2200,350){$\bullet$}

 \end{picture}
\end{center}
by the names of the nodes of the positive zoom complex $T_2$
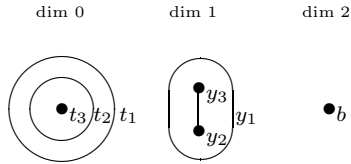
\begin{center} \xext=750 \yext=550
\begin{picture}(\xext,\yext)(\xoff,\yoff)
 \put(60,470){\tiny dim 0}
  \put(520,470){\tiny dim 1}
  \put(980,470){\tiny dim 2}
 \put(150,150){\circle{220}}
 \put(150,150){\circle{350}}
 \put(125,125){$\bullet _{t_3}$}
 \put(260,125){$ _{t_2}$}
 \put(350,125){$ _{t_1}$}

 \put(600,50){$\bullet _{y_2}$}
 \put(600,200){$\bullet _{y_3}$}
 \put(620,70){\line(0,1){150}}
 \put(630,150){\oval(220,350)}
 \put(750,125){$_{y_1}$}

 \put(1050,125){$\bullet _b$}

 \end{picture}
\end{center}
that are sent to those faces.

This duality could be compared to a restricted version of the duality between simple categories and discs \cite{Berger}, \cite{Ou}, \cite{MZ}. In that duality we have on one side some pasting diagrams described in terms of (simple) $\o$-categories and all $\o$-functors, and on the other some combinatorial structures called (finite) discs and some natural morphisms of disks. Roughly speaking, a finite disc is a finite planar tree extended by dummy/sink nodes at the ends of any linearly ordered set of sons of each node. These sink nodes do not bring any new information about the object but they are essential to get the right notion of a morphism between `such structures', i.e. all those that correspond to $\o$-functors in the dual category. If we were to throw away the sink nodes, i.e. we would consider trees instead of discs, we could still have a duality but we would need to revise the notion of a morphism on both sides. In fact, we would have duality for degeneracy maps only. On the planar tree side we would not be able to dump a `true' node onto a dummy node. This corresponds on the side of simple categories to the fact that we consider only some $\o$-functors. If we think about $\o$-functors between simple $\o$-categories as a kind of `partial composition of a {\em part} of the pasting diagram', we would need to restrict to those $\o$-functors that represent partial composition of the {\em whole} pasting diagram. In other words, if we have trees and do not have sink nodes around, our operations cannot drop any part of the pasting diagram before they start to compose them. The duality presented in this paper can be understood through this analogy. Namely, at the level of objects positive opetopes correspond to positive zoom complexes but when we look at the morphism, the natural morphisms of positive zoom complexes correspond only to degeneracies (contraction epimorphisms) on the side of positive opetopes and these maps goes in the opposite direction. This leaves of course an open question of whether we can extend positive zoom complexes one way or the other, introducing some kind of `sink nodes', so that we could have duality for more maps, e.g. all $\iota$-maps or even all $\o$-functors and not only $\iota$-epimorphisms?

The paper is organized as follows. In Section \ref{sec-tree-complexes} we define a simplification of both constellations and zoom complex that were originally introduced in \cite{KJBM}, and the maps of both constellations and positive zoom complexes. In Section \ref{sec-duality}, we describe duality for the category $\pOpeie$ of opetopes with contraction epimorphisms on one side and the category $\pZoom$ of positive zoom complexes and zoom complex maps. In Section \ref{sec-wide-complex}, we present the extension of this duality to larger categories of $\pOpeCard_\iota$ of opetopic cardinals with contraction epimorphisms and $\wZoom$ of wide positive complexes and zoom complex maps. The paper ends with an appendix where the relevant notions and facts concerning positive opetopes are recalled from \cite{Z1} and \cite{Z4}.

\section{The category of positive zoom complexes}\label{sec-tree-complexes}

\subsection{Trees}\label{sec-Trees}
A {\em tree} is a finite poset with binary sups and no infs of non-linearly ordered non-empty subsets. In particular, tree can be empty but if it is not, it has the largest element, called {\em root}.  A tree embedding is a one-to-one function that preserves and reflects order. 
We will also consider other kinds of (monotone) morphisms of trees: sup-morphisms (= preserving suprema), monotone maps (automatically preserving infs), onto maps.

{\bf Construction.}
Let $S$, $T$ be trees. Let $S_\perp$ denotes the poset obtained by adding bottom element to $S$. Then $S_\perp$ has both sups and infs, i.e. it is a lattice. If $D \subset S$ is a downward closed proper subset of $S$, then $S-D$ is again a tree. Any monotone map $f: S-D \ra T$ can be extended to an infs preserving map $f_* : S_\perp \ra T_\perp$ (sending $D$ and $\perp$ to $\perp$). $f_*$ has a left adjoint
$f^*:  T_\perp \ra S_\perp$.  If $f$ is onto, then $f^*$ reflects $\perp$, and hence it restricts to a sup-preserving morphism again named $f^*: T\ra S$.
\vskip 5mm

{\bf Some notions and notation concerning trees.}  Let $T$ be a tree,   $t,t'\in T$.
\begin{enumerate}
\item $t\prec t'$ means that $t$ is a son of $t'$ or $t'$ is the successor of $t$.

\item The suprema (infima) of a subset $X$ of a poset $T$ will be denoted by $\sup^T(X)$ or $\sup(X)$, ($\inf^T(X)$ or $\inf(X)$)  and if $X=\{ t,t'\}$, we can also write $t\vee^T t'$ of $t\vee t'$ ($t\wedge^T t'$ of $t\wedge t'$).
\item A subposet $X$ of a tree $T$ is a {\em convex subtree} of $T$  iff $X$ has the largest element, and whenever $x,x'\in X$, $s\in T$ and $x<s<x'$, then $s\in X$.
Clearly a convex subset of a tree is in particular a non-empty tree. Let $\St(T)$ denote the poset of the convex sub-trees of the tree $T$.
\item $t \perp t'$ means that either $t\leq t'$ or $t'\leq t$.
  \item $t$ is a leaf in $T$ iff the set $\{ s\in T\,:\, s\prec t\}$ is empty. $\lvs(T)$ denotes the set of leaves of the tree $T$.
  \item $\lvs^{T}(t)$ is the set of leaves of the tree $T$ over the element $t$, i.e.
\[   \lvs^{T}(t)= \{ s\in \lvs(T)\,:\, s\leq t\}. \]
  \item $\cvr^{T}(t)$ is {\em the cover of the element $t$ in the tree $T$}
  \[   \cvr^{T}(t)= \{ s\in T\,:\, s\prec t\}, \]
  i.e. the set of elements of the tree $T$ whose successor is $t$.
  \item If $X$ is a convex subtree of $T$, then  {\em the cover of $X$ in the tree $T$} is the set
  \[   \cvr^{T}(X)= \bigcup_{x\in X}  \cvr^{T}(x) - X. \]
  Note that $\cvr^{T}(t)=\cvr^{T}(\{t\})$, so the notation for value of $\cvr^{T}$ on elements and convex subsets is compatible.
\end{enumerate}

We have an easy Lemma establishing some relation between the above notions. It will be needed for the proof of the duality.

\begin{lemma}\label{lem-lvs-cov-sup}
Let $S$ be a tree.

\begin{enumerate}
  \item Let $T$ be a convex subtree of $S$ not containing leaves. Then the family set $\{ \lvs^S(s) \}_{s\in \cvr^S(T)}$ is a partition of the set $\lvs^S(\sup(T))$.
  \item Let $T$ be a convex subtree of $S$ not containing leaves and let $\{ T_i\}_{i\in I}$ be a partition of $T$. Then, for $t \in S$
  \begin{enumerate}
    \item $t\in \cvr^S(T_i)$, for some $i\in I$, iff either $t\in cvr^S(T)$ or there is ${j\in I}$ such that $t=\sup(T_j)$;
    \item $t= \sup^S(T_i)$, for some $i\in I$, iff either $t = sup(T)$ or $t\in T$ and there is $j\in I$ such that $t\in \cvr^S(T_j)$. $~~\Box$
  \end{enumerate}
\end{enumerate}
\end{lemma}

\subsection{Constellations}
 A {\em constellation} is a triple $(S_1,\sigma,S_0)$, where $S_0$ and $S_1$ are trees and $\sigma$ is a monotone function
\[ \sigma: S_1\ra \St(S_0) \]
such that
\begin{enumerate}
  \item it preserves top element;
  \item if $s,s'\in S_1$ and $\sigma(s)\cap \sigma(s')\neq\emptyset$, then $s \perp s'$.
\end{enumerate}

Let $\sigma: S_1\ra \St(S_0)$ be a constellation. Then {\em the constellation tree}  $(S_1\lhd_\sigma S_0,\leq^{co})$ is the tree arising by extension of the tree $S_1$ by nodes of the tree $S_0$ added as leaves, so that if $s_0\in S_0$ and $s_1\in S_1$, then $s_0$ is a leaf over $s_1$ iff $s_0\in \sigma(s_1)$.

Formally, the set $S_1\lhd_\sigma S_0$ is a disjoint sum of $S_1$ and $S_0$. If $s_0\in S_0$, and $s_1\in S_1$, then
the corresponding elements in $S_1\lhd_\sigma S_0$ are denoted by $s_0^\bullet$  and $s_1^\circ$, respectively.
{\em The constellation order $\leq^{co}$ in} $S_1\lhd_\sigma S_0$ is defined as follows. If $s_0,t_0\in S_0$ and $s_1,t_1\in S_1$, then
\begin{enumerate}
  \item $s_1^\circ\leq^{co} t_1^\circ$ iff $s_1\leq^{T} t_1$;
  \item $s_0^\bullet \leq^{co} s_1^\circ$ iff $s_0\in \sigma(s_1)$;
  \item $s_0^\bullet\leq^{co}t_0^\bullet$ iff $s_0=t_0$;
  \item $s_1^\circ \leq^{co}s_0^\bullet$ never holds.
\end{enumerate}
Clearly $(S_1\lhd_\sigma S_0,\leq^{co})$ is again a tree. We often drop the index $\sigma$ in $S_1\lhd_\sigma S_0$ when it does not lead to a confusion.

Let $(S_1,\sigma,S_0)$ and $(T_1,\tau,T_0)$ be two constellations.
   Any pair of tree embeddings $f_0:S_0\ra T_0$, $f_1:S_1\ra T_1$, such that $\vec{f}_0(\sigma(s))\subseteq \tau(f_1(s))$ for $s\in S_1$, induces a tree embedding of constellation orders
   $$f_1 \lhd f_0: S_1\lhd_\sigma S_0\lra T_1\lhd_\tau T_0.$$
   Such a pair is {\em a morphism of constellations}
\[ (f_1,f_0) : (S_1,\sigma,S_0) \lra (T_1,\tau,T_0)\]
 iff the induced map  $f_1 \lhd f_0$ preserves binary sups.

\vskip 2mm

{\bf Remarks and notation.}
\begin{enumerate}
\item The fibers of any constellation $\sigma$ are linearly ordered.
\item One can look at a single constellation $(S_1,\sigma,S_0)$ as data for gluing elements of a poset $S_0$ as new leaves in the posted $S_1$ along the function $\sigma$. Thus the order of $S_0$ is not essential for building one constellation order. The constellation order can be graphicly drawn with {\em leaves} from $S_0$, also called {\em vertices}, marked as dots, and inner nodes from $S_1$, also called {\em circles}, marked as circles, enclosing all the leaves under them and all the smaller circles. General elements of constellation orders are often called {\em nodes}.

  \item  Let $(T_1,\tau_0,T_0)$ and $(T_2,\tau_1,T_1)$ be two (consecutive) constellations. We can form a diagram

\begin{center} \xext=2000 \yext=150
\begin{picture}(\xext,\yext)(\xoff,\yoff)
\putmorphism(0,50)(1,0)[T_2`T_2\lhd T_1`\stackrel{_\circ}{(-)}]{600}{1}a
\putmorphism(600,50)(1,0)[\phantom{T_2\lhd T_1}`{T_1}`\stackrel{_\bullet}{(-)}]{600}{-1}a
\putmorphism(1200,50)(1,0)[\phantom{T_1}`T_1\lhd T_0`\stackrel{_\circ}{(-)}]{600}{1}a
\putmorphism(1800,50)(1,0)[\phantom{T_1\lhd T_0}`{T_0}`\stackrel{_\bullet}{(-)}]{600}{-1}a
\end{picture}
\end{center}
consisting of two embeddings of circles $(-)^\circ$ and vertices $(-)^\bullet$. Only the embeddings of circles preserve (and reflect) order.
\begin{enumerate}
  \item If $t\in T_1$, then $t^\bullet \in T_1\lhd T_0$ and  $t^\circ \in T_2\lhd T_1$. So the node $t$, depending on the order in which we consider it, can be either a vertex (leaf) or a circle (inner node). With a slight abuse we also assume that $t^\bullet =t^{\circ\bullet}\in T_1\lhd T_0$ and $t^\circ=t^{\bullet\circ}\in T_2\lhd T_1$.   We also   As we will deal with this situation very often, we will usually be careful to distinguish these two roles, when it may cause confusions, by putting either circle of dot over the node considered.
  \item Moreover, for $X\subset T_1$ we use the notation $X^\bullet \subseteq T_2\lhd T_1$ and $X^\circ=X^{\bullet\circ} \subseteq T_1\lhd T_0$.
\end{enumerate}
\end{enumerate}

\subsection{Positive zoom complexes}
{\em A positive zoom complex} $(T,\tau)$ is a sequence of constellations:
\[ \tau_0: T_1 \ra \St(T_{0}), \]
\[ \tau_1: T_2 \ra \St(T_1), \]
\[ \ldots \]
\[ \tau_{i}: T_{i+1} \ra \St(T_i), \]
\[ \ldots \]
for $i\in\o$, with almost all sets $T_i$ empty. The dimension $(T,\tau)$ is $n$ iff  $T_n$ be the last non-empty set. We write $dim(T)$ for dimension of the positive zoom complex $(T,\tau)$. $T_0$ as well as $T_{dim(T)}$ are required to be singletons.  

A {\em morphism of positive zoom complexes} $f:(S,\sigma)\ra (T,\tau)$ is a family of tree embeddings $f_i:S_i\ra T_i$, for $i\in \o$, such that,  for $i\in \o$,
\[(f_{i+1},f_i): (S_{i+1},\sigma_{i},S_i)\lra (T_{i+1},\tau_{i},T_i)\]
is a morphism of constellations, i.e. the tree embeddings
\[ \vec{f}_i=f_{i+1} \lhd f_i: S_{i+1}\lhd_{\sigma_{i} } S_i\lra T_{i+1}\lhd_{\tau_{i}} T_i \]
that preserve binary sups.

The category of positive zoom complexes and their morphisms will be denoted by $\cTree$.

\section{Duality}\label{sec-duality}

\subsection{From positive opetopes to positive zoom complexes}\label{from p-ope to t-compl}
For the notation and notions concerning positive opetopes the reader should consult Appendix and/or \cite{Z1}. In this section we define a functor
\begin{center} \xext=1200 \yext=250
\begin{picture}(\xext,\yext)(\xoff,\yoff)
\putmorphism(0,50)(1,0)[\pOpe_{\iota,epi}`\pZoom^{op}`]{800}{0}a
\putmorphism(0,50)(1,0)[\phantom{\pOpe_{\iota,epi}}`\phantom{\pZoom^{op}}`(-)^*]{800}{1}a
 \end{picture}
\end{center}

Let $P$ be a positive opetope. We shall define a positive zoom complex $(P^*,\pi)$. For $i\in\o$, the set
\[ P^*_{i}=(P_i-\gamma(P_{i+1}),\leq^-), \]
is  the set of faces of the $i$-th tree of the positive zoom complex $(P^*,\pi)$.
The $i$-th constellation map
\[ \pi_{i} : P^*_{i+1} \lra \St(P^*_i) \]
is given, for $p\in P^*_{i+1}$, by
\[ \pi_i(p) =\{ s\in P^*_i : s<^+\gamma(p) \}. \]

NB. If $p,p'\in P^*_{i+1}$, then $p<^-p'$ iff $\gamma(p)<^+\gamma(p')$.

With the notation as above, we have

\begin{proposition}
For $i\in\o$, the triple  $(P^*_{i+1}, \pi_{i+1},P^*_{i})$ defined above is a constellation.
Thus $(P^*,\pi)$ is a positive zoom complex.
\end{proposition}

{\it Proof.}~
First we show that, for $p\in P_{i+1}^*$, $\pi_i(p)$ is a convex subtree of $(P^*_{i},<^-)$. The $<^+$-least element in the $\gamma$-pencil of $\gamma\gamma(p)$ is the largest element of $\pi_i(p)$. Let $p_1,p_3\in \pi_i(p)$ and $p_2\in P^*_{i+1}$ such that $p_1<^-p_2<^-p_3$. Thus there is a maximal lower $P^*_i=P_i-\gamma(P^*_{i+1})$-path containing $p_1,p_2,p_3$. Thus, by Path Lemma (cf. \cite{Z1} or Appendix), $p_2<^+\gamma(p)$.

Next we show that $\pi_i: P_{i+1}^*\ra \St(P_i^*)$ is monotone. Let $p,p'\in P_{i+1}^*$ so that $p<^-p'$.
Then , by Proposition 5.10 of \cite{Z1}, $\gamma(p)<^+\gamma^+(p')$ and hence
\[ \pi(p) =\{ s\in P^*_i : s<^+\gamma(p) \} \subseteq \{ s\in P^*_i : s<^+\gamma(p')\}= \pi(p'), \]
as required.

Finally, we will show that if $p,p'\in P_{i+1}^*$  $s\in P^*_i$ and  $s\in \pi_i(p)\cap \pi_i(p')$, then $p\perp^- p'$. So assume that $s<^+\gamma(p)$ and $s<^+\gamma(p')$.
Let $r_1,\ldots, r_k$ be a maximal lower $P^*_{i+1}$-path such that, for some $j$, $s\in \delta(r_j)$. Then, by Lemma 5.13 of \cite{Z1},  both $p$ and $p'$ must occur in this path. So $p\perp^-p'$, as required.
$~~\Box$
\vskip 2mm

We define a poset morphism
\[ \varepsilon_{P,i} : (P^*_{i+1}\lhd P^*_i,<^{co}) \lra (P_i,<^+) \]
, for $i\in\o$, as follows. For $p\in P^*_{i+1}\lhd P^*_i$, we put
 \[  \varepsilon_{P,i}(p) \;= \; \left\{ \begin{array}{ll}
		s   & \mbox{ if } p=s^\bullet  \mbox{ for some } s\in P^*_{i}, \\
        \gamma(s) & \mbox{ if } p=s^\circ \mbox{ for some } s\in P^*_{i+1}. \\
                                    \end{array}
			    \right. \]

\begin{lemma} \label{lem-eps-adjoint}
The morphism $\varepsilon_{P,i}$ defined above is an order isomorphism, for $i\in\o$.
\end{lemma}

{\it Proof.}~
If $p_1,p_2\in P^*_{i+1}=P_{i+1}-\gamma(P_{i+2})$ and $p_1\neq p_2$, then $\gamma(p_1)\neq\gamma(p_2)$.
Moreover, if $p\in  P^*_{i}=P_{i}-\gamma(P_{i+1})$, then $\varepsilon_{P,i}(p^\bullet)=p\neq \gamma(p_1)= \varepsilon_{P,i}(p_1^\circ)$. Thus $\varepsilon_{P,i}$ is one-to-one. It is onto as well, since, by Proposition 5.19.3 of \cite{Z1}, $P_i=(P_i-\gamma(P_{i+1}) \cup (\gamma(P_{i+1}-\gamma(P_{i+2}))$.

It remains to show that $\varepsilon_{P,i}$ preserves and reflects order. We have that
\[ p^\bullet\; <^{co}\; p_1^\circ\;\; {\rm iff}\;\; p\;<^+\; \gamma(p_1)\;\; {\rm iff}\;\; \varepsilon_{P,i}(p^\bullet)\;<^+\; \varepsilon_{P,i}(p_1^\circ). \]
Moreover, using Lemma 5.9.6 of \cite{Z1}, we have
\[ p_1^\circ \; <^{co}\; p_2^\circ\;\; {\rm iff}\;\; p_1\;<^-\; p_2\;\; {\rm iff}\;\; \gamma( p_1)\;<^+\; \gamma(p_2)\;\; {\rm iff}\;\; \varepsilon_{P,i}(p_1^\circ)\;<^+\; \varepsilon_{P,i}(p_2^\circ). \]
The other two cases are obvious.
$~~\Box$

Let $f: P\ra Q$ be a $\iota$-epimorphism of positive opetopes. We define a map of positive zoom complexes
\[ f^*=\{f^*_i\}_{i\in \o} : (Q^*,\pi)\lra (P^*,\pi). \]
For $i\in \o$, the map $f^*_i : Q^*_i\ra P^*_i$ is defined as follows.
 Let $q\in Q^*_i=Q_i-\gamma(Q_{i+1})$,  $p\in P^*_i=P_i-\gamma(P_{i+1})$, $0\leq i$. Then
\[ f^*_i(q) = p\]
iff $p$ is the unique element of $P_i-\gamma(P_{i+1})$ so that $f(p)=q$. Such an element exists since $f_i$ is epi.

We can also describe the above map using the construction from Section \ref{sec-Trees}, as follows. We have a monotone onto map $f_i : (P_i -\ker(f),\leq^+) \ra (Q_i,\leq^+)$. As $P_i\cap \ker(f)$ is a proper downward closed subset of $P_i$, $f_i$ extends to an all infs preserving map $f_{i,*}: P_{i,\perp}\ra  Q_{i,\perp}$ sending $P_i \cap \ker(f)$ to $\perp$. Thus it has a left adjoint $\bar{f}_{i}: Q_{i,\perp}\ra  P_{i,\perp}$. Clearly,  $\bar{f}_{i}$, defined this way,   preserves sups. For $q\in Q_i$, $\bar{f}_{i}(q)$ picks the $<^+$-least element in the fiber of the function $f_i$ over element $q$. We have

\begin{lemma} \label{lem-eps}
With the notation as above, for $i\in \o$, the following diagram
\begin{center} \xext=1200 \yext=700
\begin{picture}(\xext,\yext)(\xoff,\yoff)
\setsqparms[-1`1`1`-1;1400`500]
\putsquare(0,50)[(P^*_{i+1}\lhd P^*_i)_\perp`(Q^*_{i+1}\lhd Q^*_i)_\perp`(P_i)_\perp`(Q_i)_\perp;(f^*_{i+1}\lhd f^*_i)_\perp`(\varepsilon_{P,i})_\perp`(\varepsilon_{Q,i})_\perp`\bar{f}_{i}]
\end{picture}
\end{center}
is well defined and commutes. In particular, $f^*_{i+1}\lhd f^*_i$ preserves binary sups.
\end{lemma}

{\it Proof.}~ First we shall verify that, for $i\in\o$, $(f^*_{i+1},f^*_i)$ induce a monotone function
\begin{center} \xext=1400 \yext=100
\begin{picture}(\xext,\yext)(\xoff,\yoff)
 \putmorphism(0,50)(1,0)[(P^*_{i+1}\lhd P^*_i)`(Q^*_{i+1}\lhd Q^*_i)`(f^*_{i+1}\lhd f^*_i)]{1400}{-1}a
\end{picture}
\end{center}
that is, for $q\in Q^*_{i+1}$
\[ \vec{f}^*_i(\pi_i(q)) \subseteq \pi_i(f^*_{i+1}(q)). \]
Let $q'\in \pi(q)\subseteq Q^*_i$  and let $f^*_{i+1}(q)=p\in P^*_{i+1}$ and $f^*_{i}(q')=p'\in P^*_{i}$. Then $q'<^+\gamma(q)$, $f_{i+1}(p)=q$ and $f_{i}(p')=q'$. We have
\[ p'=f_i(q')\leq^+ f_i(\gamma(q)) = \gamma(f_{i+1}(q)) = \gamma(p), \]
that is $p'\in \pi(p)$, as required.

Next we show that the square commutes. Let  $q\in Q^*_{i}$, $p\in P^*_{i}$ so that $q^\bullet\in Q^*_{i+1}\lhd Q^*_i$ and let $f^*_{i+1}\lhd f^*_i(q^\bullet) = p^\bullet\in P^*_{i+1}\lhd P^*_i$. Thus $f_{i,*}(p)=q\in Q^*_{i}$ and, as $p$ is a leaf,  $\bar{f}_i(q)=p$. Hence the square commutes in this case.

Now, let $q\in Q^*_{i+1}$, $p\in P^*_{i+1}$, so that  $q^\circ\in Q^*_{i+1}\lhd Q^*_i$ and  $f^*_{i+1}\lhd f^*_i(q^\circ) = p^\circ\in P^*_{i+1}\lhd P^*_i$. Thus $f_{i+1,*}(p)=q$. So $p$ is indeed in the fiber of $f_{i+1,*}$ over $q$  and $p\not\in ker(f)$. We need to show that $\gamma(p)$ is $<^+$-minimal in the fiber of $f_i$ over $\gamma(q)$. Suppose to the contrary that there is $p'\in  P_i$ such that $p'<^+ \gamma(p)$ and yet $f_{i,*}(p')=\gamma(q)$. Then, as $p\in P_{i+1}-\gamma(P_{i+2})$ and $p'<^+ \gamma(p)$, there is a $p''\in \delta(p)$ such that $p'\leq^+p''$. We have $$\gamma(q)=f_{i,*}(p')\leq^+f_{i,*}(p'')\leq^+f_{i,*}(\gamma(p))=\gamma(q).$$
But then $f_{i+1,*}(p'')= q = f_{i+1,*}(\gamma(p))$ and hence $p\in \ker(f)$, contrary to the supposition. Thus the diagram commutes in this case, as well.

Finally, since $\bar{f}_{i}$ is a left adjoint it preserves sups, and, by the above, $(f^*_{i+1}\lhd f^*_i)_\perp$ preserves sups as well. But then $(f^*_{i+1}\lhd f^*_i)$ preserves non-empty sups, as required.  $~~\Box$
\vskip 2mm

As a corollary of Lemma \ref{lem-eps}, we get

\begin{proposition}
Let $f: (P,\gamma,\delta)\ra (Q,\gamma,\delta)$ be an epi $\iota$-map of positive opetopes. Then the family of maps $f^*=\{f^*_i\}_{i\in \o} : (Q^*,\pi)\ra (P^*,\pi)$ defined above is a morphism of positive zoom complexes.
\end{proposition}

{\it Proof.}~ From Lemma \ref{lem-eps} follows that
\[  (f_{i+1},f_i): (Q^*_{i+1},\pi_i,Q^*_i)\lra (P^*_{i+1},\pi_i,P^*_i) \]
is a constellation morphism for and $i\in \o$.
$~~\Box$
\vskip 2mm

{\bf Examples.} We explain below in more detail the correspondence sketched in the introduction between the $\iota$-epimorphism $f:O_3\ra O_2$ and its dual positive zoom complex embedding.
\begin{enumerate}
  \item The dual of the opetope $Q=O_2$ of dimension 2
  \begin{center} \xext=500 \yext=400
\begin{picture}(\xext,\yext)(\xoff,\yoff)
 \settriparms[-1`1`1;300]
 \putAtriangle(50,50)[t_2`t_3`t_1;y_3`y_2`y_1]
 \put(300,150){$\Da\! b$}
 \end{picture}
\end{center}
is the positive zoom complex $Q^*=T_2$ with nodes
\[ Q^*_0 =\{ t_3 \}, \]
\[ Q^*_1 , =\{ y_2> y_3\}, \]
\[ Q^*_2 , =\{ b \}, \]
and the constellation maps
\[ \pi_0(y_2)=\pi_0(y_3)=\{ t_3 \}, \]

\[ \pi_1(b)=\{ y_2, y_3 \}. \]
Such a positive zoom complex $Q^*$ can be drawn as follows:

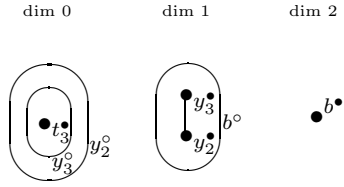
\begin{center} \xext=750 \yext=550
\begin{picture}(\xext,\yext)(\xoff,\yoff)
 \put(60,500){\tiny dim 0}
  \put(540,500){\tiny dim 1}
  \put(980,500){\tiny dim 2}
 \put(150,130){\oval(150,240)}
 \put(150,130){\oval(270,400)}
 \put(110,100){$\bullet _{t_3^\bullet}$}
 \put(160,-7){$ _{y_3^\circ}$}
 \put(290,50){$ _{y_2^\circ}$}

 \put(600,60){$\bullet _{y_2^\bullet}$}
 \put(600,200){$\bullet _{y_3^\bullet}$}
 \put(620,70){\line(0,1){150}}
 \put(630,150){\oval(220,370)}
 \put(750,125){$_{b^\circ}$}

 \put(1050,125){$\bullet ^{b^\bullet}$}
 \end{picture}
\end{center}

  \item The dual of the opetope $P=O_3$ of dimension 3
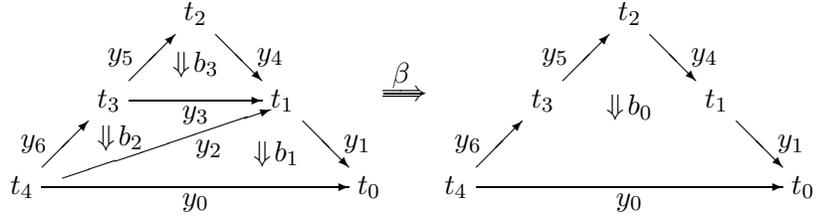
\begin{center} \xext=2400 \yext=800
\begin{picture}(\xext,\yext)(\xoff,\yoff)

 \settriparms[-1`1`1;300]
 \putAtriangle(300,350)[t_2`t_3`t_1;y_5`y_4`y_3]
 \put(520,450){$\Da\! b_3$}
 \settriparms[-1`0`0;300]
 \putAtriangle(0,50)[\phantom{t_3}`t_4`;y_6``]
 \settriparms[0`1`0;300]
 \putAtriangle(600,50)[\phantom{t_1}``t_0;`y_1`]
 \putmorphism(0,50)(1,0)[\phantom{t_4}`\phantom{t_0}`y_0]{1200}{1}b
 \putmorphism(350,150)(3,1)[\phantom{t_4}`\phantom{t_0}`]{300}{1}a
 \put(600,170){$y_2$}
 \put(260,195){$\Da\! b_2$}
 \put(800,140){$\Da\! b_1$}

\put(1240,350){$\Longrightarrow$}
  \put(1250,375){\line(1,0){135}}
   \put(1280,410){$\beta$}

\settriparms[-1`1`0;300]
 \putAtriangle(1800,350)[t_2`t_3`t_1;y_5`y_4`]
 \put(2020,300){$\Da\! b_0$}
 \settriparms[-1`0`0;300]
 \putAtriangle(1500,50)[\phantom{t_3}`t_4`;y_6``]
 \settriparms[0`1`0;300]
 \putAtriangle(2100,50)[\phantom{t_1}``t_0;`y_1`]
 \putmorphism(1500,50)(1,0)[\phantom{t_4}`\phantom{t_0}`y_0]{1200}{1}b
\end{picture}
\end{center}
is a positive zoom complex $P^*$ with nodes
\[ P^*_0 =\{ t_4 \}, \]
\[ P^*_1 , =\{ y_1> y_4 > y_5 > y_6 \}, \]
\[ P^*_2 , =\{ b_1 > b_2 > b_3 \}, \]
\[ P^*_3 , =\{ \beta\}, \]
and the constellation maps
\[ \pi_0(y_1)=\pi_0(y_4)=\pi_0(y_5)=\pi_0(y_6)=\{ t_4 \}, \]

\[ \pi_1(b_1)=\{ y_1, y_4, y_5, y_6, \}, \]
\[ \pi_1(b_2)=\{  y_4, y_5, y_6, \}, \]
\[ \pi_1(b_3)=\{ y_4, y_5, \}, \]

\[ \pi_2(\beta)=\{ b_1, b_2, b_3. \}. \]
Such a positive zoom complex $P^*$ can be drawn as follows:

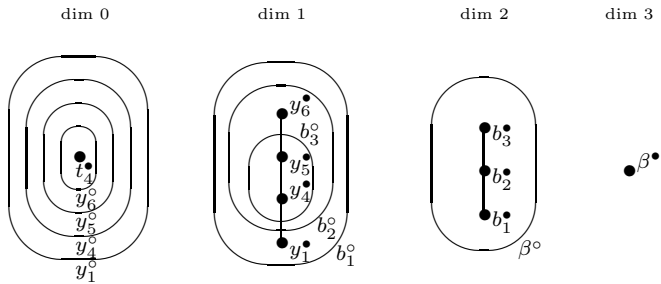
\begin{center} \xext=2200 \yext=950
\begin{picture}(\xext,\yext)(\xoff,\yoff)
 \put(260,900){\tiny dim 0}
  \put(940,900){\tiny dim 1}
  \put(1640,900){\tiny dim 2}
  \put(2140,900){\tiny dim 3}
 \put(300,400){$\bullet$}
 \put(310,370){$ _{t_4^\bullet}$}
  \put(320,420){\oval(120,220)}
  \put(310,280){$ _{y_6^\circ}$}
  \put(320,420){\oval(240,380)}
  \put(310,198){$ _{y_5^\circ}$}
   \put(320,420){\oval(360,540)}
  \put(310,115){$ _{y_4^\circ}$}
   \put(320,420){\oval(480,700)}
  \put(310,30){$ _{y_1^\circ}$}

  \put(1000,100){$\bullet _{y_1^\bullet}$}
 \put(1000,255){$\bullet ^{y_4^\bullet}$}
  \put(1000,400){$\bullet _{y_5^\bullet}$}
 \put(1000,550){$\bullet ^{y_6^\bullet}$}
 \put(1020,120){\line(0,1){450}}

 \put(1020,350){\oval(220,300)}
 \put(1085,505){$ _{b_3^\circ}$}

 \put(1020,420){\oval(360,500)}
 \put(1145,180){$ _{b_2^\circ}$}

  \put(1020,400){\oval(460,700)}
 \put(1210,85){$ _{b_1^\circ}$}

  \put(1700,200){$\bullet _{b_1^\bullet}$}
 \put(1700,350){$\bullet _{b_2^\bullet}$}
  \put(1700,500){$\bullet _{b_3^\bullet}$}
 \put(1720,220){\line(0,1){300}}

   \put(1720,400){\oval(360,600)}
 \put(1840,100){$ _{\beta^\circ}$}

 \put(2200,350){$\bullet ^{\beta^\bullet}$}

 \end{picture}
\end{center}

  \item The dual of the  $\iota$-epimorphism $f:P\ra Q$ from opetope $P=O_3$ to the opetope $Q=O_2$ given by
  \[ f_0(t_0)=f_1(y_1)=f_0(t_1)=t_1,  \]
  \[ f_0(t_2)=f_1(y_5)=f_0(t_3)=t_2,\]
  \[ f_0(t_4)=t_3, \]

  \[  f_1(y_2)=f_2(b_1)=f_1(y_0)=y_1,  \]
  \[ f_1(y_4)=f_2(b_3)=f_1(y_3)=y_2, \]
  \[ f_1(y_6)=y_3, \]

  \[ f_2(b_2)=f_3(\beta)=f_2(b_0)=b, \]
  is a morphism of positive zoom complexes $f^*:Q^*\ra P^*$ such that
  \[ f^*_0(t_3)=t_4,  \]

  \[ f^*_1(y_2)=y_4, \]
  \[ f^*_1(y_3)=y_6,\]

  \[ f^*_2(b)=b_2.\]
\end{enumerate}

\subsection{From  positive zoom complexes to positive opetopes}

In this section we define a functor 
\begin{center} \xext=1200 \yext=250
\begin{picture}(\xext,\yext)(\xoff,\yoff)
\putmorphism(0,50)(1,0)[\pOpe_{\iota,epi}`\pZoom^{op}`]{800}{0}a
\putmorphism(0,50)(1,0)[\phantom{\pOpe_{\iota,epi}}`\phantom{\pZoom^{op}}`(-)^*]{800}{-1}b
 \end{picture}
\end{center}

Let $(S,\sigma)=\{S_i,\sigma_i\}_{i\in\o}$ be a positive zoom complex. We define the positive opetope $\{ S^*_i,\gamma^i,\delta^i\}_{i\in\o}$, as follows.
We put
\[ (S^*_i,<^{co})=(S_{i+1}\lhd_{\sigma_i}S_{i},<^{co}), \]
i.e. the set $S^*_i$  of $i$-dimensional faces of the positive opetope $S^*$ is the universe of the $i$-th constellation poset of $(S,\sigma)$. Later we shall prove that the constellation order $<^{co}$ agree with the upper order $<^+$, defined using the operations $\gamma$ and $\delta$ below.

Let $i\in\o$. The $i-${\em th codomain operation}

\[ \gamma :  S^*_{i+1} \lra S^*_{i}\]
is defined, for $p\in S^*_{i+1}$, as follows
\[ \gamma(p) \;= \;  \sup\!\,^{S^*_{i}}(\lvs^{S^*_{i+1}}(p)^\circ). \]
In words
\begin{enumerate}
  \item if the face $p$ is a vertex, i.e. $p=t^\bullet$ for some $t\in S_i$, then its codomain $\gamma(t^\bullet)=t^\circ$, i.e. it is `the same' $t$ but considered as a circle one dimension below;
  \item if the face $p$ is a circle for some $t\in S_{i+1}$, i.e. $p=t^\circ$, then its codomain $\gamma(p)$ is the circle $s^\circ$ whose corresponding node $s^\bullet$ is the supremum $\sup^{S^*_i}(\lvs^{S^*_{i+1}}(t^\circ)^\circ)$ in $S^*_{i}$ of the leaves in $S^*_{i+1}$ over $t^\circ$ considered as circles one dimension below, in $S^*_{i}$.
\end{enumerate}

The $i-${\em th domain operation}
\[ \delta :  S^*_{i+1} \lra {\cP}_{\neq\emptyset}(S^*_{i}) \]
is defined, for $p\in S^*_{i+1}$, as follows
\[ \delta(p) \;= \; \cvr^{S^*_{i}}(\lvs^{S^*_{i+1}}(p)^\circ).\]
In words
\begin{enumerate}
  \item if the face $p$ is a vertex, i.e. $p=t^\bullet$ for some $t\in S_i$, then its domain $\delta(t^\bullet)=\cvr^{S^*_i}(t^\circ)$, i.e. it is the cover of `the same' $t$ but considered as a circle one dimension below;
  \item if the face $p$ is a circle, i.e. $p=t^\circ$ for some $t\in S_{i+1}$, then its domain $\delta(t^\circ)$ is the sum of $\delta$'s applied to the leaves/vertices over $t^\circ$ in $S^*_{i+1}$ considered as circles one dimension below in $S^*_i$ minus these leaves considered as circles.
\end{enumerate}

\begin{lemma}\label{opetope-from-ctree-ob}
  Let $(S,\sigma)$ be a positive zoom complex. Then the face structure $(S,\sigma)^*=(S^*,\gamma,\delta)$ defined above is a positive opetope.
\end{lemma}

{\it Proof.}~
Let $(S^*,\gamma,\delta)$ a face structure as defined above. We shall check that it satisfies the axioms of positive opetopes.
\vskip 2mm

{\bf Globularity.} We shall use Lemma \ref{lem-lvs-cov-sup}.

Fix $s\in S^*_{i+2}$, for some $i\geq 0$. Then $\lvs^{i+2}(s)^{\circ}$ is a convex subtree of  $S^*_{i+1}$ not containing leaves. Thus, by Lemma \ref{lem-lvs-cov-sup}.1, the family of sets
$$\{ \lvs^{i+1}(r)\}_{r\in \cvr^{i+1}(\lvs^{i+2}(s)^{\circ})}$$
is a partition of the set
$$\lvs^{i+1}(\sup\!^{i+1}(\lvs^{i+2}(s)^{\circ})).$$

If $r=r'^{\bullet}\in \cvr^{i+1}(\lvs^{i+2}(s)^{\circ})$, for some $r'\in S_{i+1}$, then  $\lvs^{i+1}(r)^\circ = r'^{\circ}$.
If $r=r'^{\circ} \in \cvr^{i+1}(\lvs^{i+2}(s)^{\circ})$, for some $r'\in S_{i+2}$, then  $\lvs^{i+1}(r)^\circ = \sigma_{i+1}(r)$. Thus in any case $\lvs^{i+1}(r)^\circ$ is a convex subtree, and hence the partition, give rise to the partition of a convex subtree
$$\lvs^{i+1}(\sup\!^{i+1}(\lvs^{i+2}(s)^{\circ}))^\circ,$$
into a family of convex subtrees
$$\{ \lvs^{i+1}(r)\}_{r\in \cvr^{i+1}(\lvs^{i+2}(s)^{\circ})}^\circ$$
of the tree $S^*_i$.

Then, using Lemma \ref{lem-lvs-cov-sup}.2, we get
\[ \gamma\gamma(s)=   \]

\[ = \sup\!^i(\lvs^{i+1}(\sup\!^{i+1}(\lvs^{i+2}(s)^\circ))^\circ)    =\]

\[ = \bigcup_{r\in \cvr^{i+1}(\lvs^{i+2}(s)^\circ)} \sup\!^{i}(\lvs^{i+1}(r)^\circ) -  \bigcup_{r\in \cvr^{i+1}(\lvs^{i+2}(s)^\circ)} \cvr^{i}(\lvs^{i+1}(r)^\circ) =\]

\[ = \gamma\delta(s)-\delta\delta(s), \]
and
\[ \delta\gamma(s)= \]

\[ = \cvr^i(\lvs^{i+1}(\sup\!^{i+1}(\lvs^{i+2}(s)^\circ))^\circ)    =\]

\[ = \bigcup_{r\in \cvr^{i+1}(\lvs^{i+2}(s)^\circ)} \cvr^{i}(\lvs^{i+1}(r)^\circ) - \bigcup_{r\in \cvr^{i+1}(\lvs^{i+2}(s)^\circ)} \sup\!^{i}(\lvs^{i+1}(r)^\circ) =\]

\[ = \delta\delta(s)-\gamma\delta(s), \]
 as required.

\vskip 2mm

{\bf Strictness.}

We shall show that the transitive relation  $<^+$, defined using $\gamma$'s and $\delta$'s, coincides with the constellation order $<^{co}$.

Let $s, s'^\circ \in S_{i+1}\lhd S_{i}$, for some $0\leq i \leq dim(S)$.
Then $\gamma(s'^\bullet)= s'^\circ$. Moreover, $s\prec^{co} s'^\circ$ iff $s\in \cvr^i(s'^\circ) = \delta(s'^\bullet)$.
The latter condition means that $s<^+ s'^\circ$. Thus $\prec^{co}\subseteq <^+$.

It remains to show that $<^+\subseteq <^{co}$. Assume $s,s'\in S_{i+1}\lhd S_i$ and that there is $r\in S_{i+2}\lhd S_{i+1}$ such that $s\in \delta(r)$ and $\gamma(r)=s'$. We shall show that $s<^{co}s'$.

If $r=r'^\bullet$, for some $r'\in S_{i+1}$, then $s'=r'^\circ$ and $s\in \cvr^i(r'^\circ)$ so $s\prec^{co} s'$ indeed.

Now assume that $r=r'^\circ$, for some $r'\in S_{i+2}$. Let $s=s_0\in \delta(r'^\circ) = \cvr^i(\lvs^{i+1}(r'^\circ)^\circ)$. Let $s_1$ be the $<^{co}$-successor of $s_0$, i.e. $s_0\prec^{co} s_1$. Since $\lvs^{i+1}(r'^\circ)$ is a convex tree, there is a path $s_1,\ldots,s_k$ in $\lvs^{i+1}(r'^\circ)^\circ$ such that $s_i\prec^{co} s_{i+1}$ and
$$s_k= \sup^i(\lvs^{i+1}(r'^\circ)^\circ) = \gamma(r^\circ)=s'.$$
Thus $s=s_0<^{co} s_k=s'$, as required.

\vskip 2mm

{\bf Disjointness.}

Let $s,t\in S^*_{i+1}$. If $s<^+t$, then $\lvs^{i+1}(s)\subseteq \lvs^{i+1}(t)$. On the other hand if $s=s_0,\ldots, s_k=t$ is a lower path in $S^*_{i+1}$, i.e.,
$$ \gamma(s_i) =   \sup\!\,^{S^*_{i}}(\lvs^{S^*_{i+1}}(s_i)^\circ) \in    \cvr^{S^*_{i}}(\lvs^{S^*_{i+1}}(s_{i+1})^\circ)=\delta(s_{i+1}),$$
for $i=0,\ldots, k-1$.  In other words, $\gamma(s_i)$, the largest element of  $\lvs^{S^*_{i+1}}(s_i)^\circ$, is smaller than the least element in $\lvs^{S^*_{i+1}}(s_{i+1})^\circ$ comparable with $\gamma(s_i)$, . Thus the elements of the sets $\{ \lvs^{S^*_{i+1}}(s_i)^\circ\}_{i=1,\ldots, k}$ are pairwise disjoint.
In particular the sets $\lvs^{S^*_{i+1}}(s)$ and $\lvs^{S^*_{i+1}}(t)$ are disjoint whenever $s\perp^-t$. This the orders $<^-$ and $<^+$ are disjoint, as required.

\vskip 2mm

{\bf Pencil linearity.}

Let $s,t \in S^*_i$, for some $0\leq i\leq dim(S)$.

Assume that $s\neq t$ and $\gamma(s)=\gamma(t)$.  Then $s$ and $t$ cannot be leaves at the same time.
If $s$ is a leaf, then $s<^{co}t$ and hence $s<^+t$, by the above.
If both $s$ and $t$ are inner nodes, 
then
\[ \sigma_i(s)^\circ\ni \sup^{i-1}(\lvs^i(s)^\circ) =  \sup^{i-1}(\lvs^i(t)^\circ)\in \sigma_{\it i}(t)^\circ. \]
Thus $\sigma_i(s)\cap\sigma_i(t)\neq \emptyset$, and, as $\sigma_i$ is a constellation, we have $s\perp^+ t$.

Now assume that there is $r\in \delta(s)\cap\delta(t)$. Let $s_1^\circ$ be the successor of $r$, i.e. $r\prec^{co} s_1^\circ$.
Hence $s_1^\circ\in \lvs^i(s)^\circ\cap\lvs^i(t)^\circ$ and hence
\[ s_1^\bullet\in \lvs^i(s)\cap\lvs^i(t). \]
If $s$ and $t$ were leaves, then we would have $s=t$.

If $s$ is a leaf and $t$ is an inner node, then $s\in \lvs^i(t)$ and hence $s<^+t$.

If both $s$ and $t$ are inner nodes, then
\[ s_1^\bullet\in \lvs^i(s)\cap \lvs^i(t) =\sigma_{i+1}(s)^\bullet\cap\sigma_{i+1}(t)^\bullet, \]
and as $\sigma_{i+1}$ is a constellation, $s\perp^+t$.
$~~\Box$

Let $f :(S,\sigma) \ra (T,\tau)$ be a map of positive zoom complexes. It gives rise to maps of faces,  for $i\in\o$,
\[   \vec{f}_i=f_{i+1}\lhd f_{i} :  S^*_i=S_{i+1}\lhd S_{i} \lra T^*_i=T_{i+1}\lhd T_{i}, \]
that, by definition, preserve binary sups.  Note that the maps $\vec{f}_i$'s do not preserve the domains or codomains just defined above, in general.
These maps induce the $\iota$-epimorphism of positive opetopes $f^*:T^*\ra S^*$, i.e., the maps $f^*_i: T^*_i\ra S^*_{\leq i}$, for $i\in \o$, as follows. Let $t\in T^*_i$ and $s\in S^*_j$, with $0\leq j\leq i \leq dim(T)$.
 Then
 \[ f^*_i(t) = s \]
 iff
 \begin{enumerate}
   \item $j$ is the maximal number such that there is $s'\in S^*_j$ that $$\vec{f}_j(s')\leq^{co}\delta\gamma^{(j+1)}(t);$$
   \item and $s$ is the $\leq^{co}$-maximal $s'\in S^*_j$ satisfying the above inequality.
 \end{enumerate}

\begin{lemma}\label{opetope-from-ctree-mor}
  Let $f:(S,\sigma)\ra (T,\tau)$ be a morphism of positive zoom complexes. Then the set of maps $f^*=\{f^*_i\}_{i\in\o} : (T,\tau)^*\ra (S,\sigma)^*$ is a $\iota$-epimorphism of positive opetopes.
\end{lemma}

{\it Proof.}~
Let us fix a morphism of positive zoom complexes $f:(S,\sigma)\ra (T,\tau)$,  $i\in \o$, $s\in S^*_i=S_{i+1}\lhd S_i$ and $t=\vec{f}_i(s)$. Then $\vec{f}_i(s)\leq^{co} t$ and, as $\vec{f}_i$ is one-to-one, it is the largest such $s$. Thus $f^*_i(t)=s$. Since $s$ was arbitrary,  $f^*_i$ is onto, for any $i\in \o$ and hence $f^*$ epi.

For  preservation of both codomains and domains by $f^*$, we fix $i>0$ and $t\in T^*_i$  and we consider three cases:
\begin{enumerate}
  \item $f^*_i(t)\in  S^*_i$;
  \item $f^*_i(t)\in  S^*_{i-1}$;
  \item $f^*_i(t)\in  S^*_{j}$, for some $j<i-1$.
\end{enumerate}

{\bf Preservation of codomains $\gamma$}.
\vskip 2mm
{\em Case $\gamma$.1}: $f^*_i(t)=s\in  S^*_i$.

First we shall show that $\vec{f}_{i-1}(\gamma(s))\leq^{co} \gamma(t)$. Since $\vec{f}_i$'s are monotone and preserve leaves, we have
\[ \vec{f}_{i}({\rm lvs^i}(s)) \subseteq {\rm lvs^i}(\vec{f}_i(s)) \subseteq \lvs^i(t).\]
Using the above and the fact that $\vec{f}_i$'s preserve sups, we have
\[ \vec{f}_{i-1}(\gamma(s)) =  \]
\[= \vec{f}_{i-1}(\sup\!^{i-1}(\lvs^i(s)^\circ)) = \]
\[ ={\rm sup}\!^{i-1}(\vec{f}_{i-1}(\lvs^i(s)^\circ)) \leq^{co}\]
\[ \leq^{co} {\rm sup}\!^{i-1}(\lvs^i(t)^\circ))= \]
\[=\gamma(t).  \]
Now, contrary to the claim we want to prove, we assume that there is  $s_1^\circ\in S^*_{i-1}$ such that $\gamma(s)<^{co}s_1$ and
\[ \vec{f}_{i-1}(\gamma(s)) <^{co} \vec{f}_{i-1}(s_1^\circ) \leq^{co} \gamma(t)\in\lvs^i(t)^\circ\subseteq T^*_{i-1}. \]
Thus  $ \vec{f}_{i}(s_1^\bullet)\in T^*_i-\vec{f}_{i}(\lvs^i(s)^\circ)$. As $\vec{f}_{i}(s)\leq^{co} t$, we have $\vec{f}_{i}(\lvs^i(s))\subseteq \lvs^i(t)$. Hence $\vec{f}_{i-1}(\gamma(s))\in \lvs^i(t)^\circ$. Since $\lvs^i(t)^\circ$ is a convex subtree, we have
\[ \vec{f}_{i-1}(s_1^\circ)\in \lvs^i(t)^\circ. \]

As $\vec{f}_i(s_1^\bullet)^\circ=\vec{f}_{i-1}(s_1^\circ)$, we have $\vec{f}_{i}(s_1^\bullet)\leq^{co} t$. Since we also have $s_1^\bullet \not\in \lvs^i(s)$,
we get that
\[ s<\sup\!^i(\{ s_1^\bullet\}\cup \lvs^i(s))=s_2^\circ \]
and
\[   \vec{f}_{i}(s) <^{co}  \vec{f}_{i}(s_2^\circ) =\]
\[ = \vec{f}_{i}(\sup\!^i(\{ s_1^\bullet\}\cup \lvs^i(s))) =\]
\[ =  {\rm sup^i}(\{\vec{f}_{i}(s_1^\bullet) \}\cup \vec{f}_{i}(\lvs^i(s)))\leq^{co} t.\]
This is a contradiction with the fact that $f^*_i(t)=s$. This ends the proof of Case $\gamma.1$.

\vskip 2mm
{\em Case $\gamma$.2}:  $f^*_i(t)=s_1\in  S^*_{i-1}$.

Thus we have a $t_1\in \delta(t)$ such that
\[ \vec{f}_{i-1}(s_1) \leq^{co} t_1 <^{co} \gamma(t). \]
We need to show that $s_1$ is the largest such an element of $S^*_{i-1}$ that $\vec{f}_{i-1}(s_1) \leq^{co} \gamma(t)$. Suppose to the contrary that there is $s_2^\circ\in S^*_{i-1}$ such that $s_1<^{co}   s_2^\circ$ and
\[  \vec{f}_{i-1}(s_1)<^{co}  \vec{f}_{i-1}(s_2^\circ) \leq^{co} \gamma(t).\]
We have $\vec{f}_{i-1}(s_1)\leq^{co} t_1$ and $\vec{f}_{i-1}(s_1)<^{co} \vec{f}_{i-1}(s_2^\circ)$, and, as we cannot have $\vec{f}_{i-1}(s_2^\circ)\leq^{co}t_1$,
we have
\[ \vec{f}_{i-1}(s_1)\leq^{co} t_1<^{co} \vec{f}_{i-1}(s_2^\circ) \leq^{co} \gamma(t).\]
Since $\lvs^i(t)^\circ$ is a convex subtree of $T^*_{i-1}$, it follows that $f^*_{i-1}(s_2^\circ)\in \lvs^i(t)^\circ$. Thus $f^*_{i}(s_2^\bullet)\in\lvs^i(t)$, i.e. $f^*_{i}(s_2^\bullet)\leq^{co} t$. Hence $f^*_i(t)\in S^*_i$, contrary to the supposition. This ends the proof of Case $\gamma.2$.

\vskip 2mm
{\em Case $\gamma$.3}: $f^*_i(t)=s_1\in  S^*_{j}$, for some $j<i-1$.

Suppose there is $s_2\in S^*_{i-1}$ such that $\vec{f}_{i-1}(s_2)\leq^{co}\gamma(t)$. Then, as $j<i-1$, $\vec{f}_{i-1}(s_2)\not\leq^{co}t'$, for all $t'\in\delta(t)$.
Thus there is $t_1\in\delta(t)$ and  $s_3\in S^*_i$ so that $s_2=s_3^\circ$
and
\[t_1<^{co}\vec{f}_{i-1}(s_3^\circ)\leq^{co}\gamma(t).\]
As $\lvs^i(t)^\circ$ is a convex subtree, we have $\vec{f}_{i-1}(s_3^\circ)\in\lvs^i(t)^\circ$ and then $\vec{f}_i(s_3^\bullet)\in \lvs^i(t)$, i.e., $\vec{f}_i(s_3^\bullet)\leq^{co} t$ . This contradicts the fact that $j<i-1$.
Thus $f^*_{i-1}(\gamma(t))=s_4\in S^*_{j'}$ such that $j'<i-1$.
If $j'>j$, then
\[ \vec{f}_{j'}(s_4)\leq^{co}\delta\gamma^{(j'+1)}(\gamma(t)) =\delta\gamma^{(j'+1)}(t)\]
and this contradicts the choice of $s_1\in S^*_j$. If $j>j'$, then
\[ \vec{f}_{j}(s_1)\leq^{co}\delta\gamma^{(j+1)}(t)=\delta\gamma^{(j+1)}(\gamma(t)) \]
contradicting the choice of $s_4\in S^*_{j'}$. Thus $j=j'$ and $s_1=s_3$, as required. This ends the proof of Case $\gamma.3$.
\vskip 3mm
{\bf Preservation of domains $\delta$}.

\vskip 2mm
{\em Case $\delta$.1}: $f^*_i(t)=s\in  S^*_i$.

We shall show that $f^*_i$ restricts to a bijection
\[ f^*_{i-1\lc t} : \delta(t)-\ker(f^*)\lra \delta(s). \]
Let $t_1\in \delta(t)-\ker(f^*)$.  Thus there is $s_1\in S^*_{i-1}$ such that $f^*_{i-1}(t_1)=s_1$ and hence $\vec{f}_{i-1}(s_1)\leq^{co}t_1$. Since $f^*$ preserves codomains $\gamma(s)=f^*_{i-1}(\gamma(t))$.

Since $\vec{f}_{i-1}(s_1)\not\in{\rm lvs^i}(t)^\circ\supseteq \vec{f}_{i-1}(\lvs^i(s)^\circ)$, it follows that $s_1\not\in \lvs^i(s)^\circ$.

We shall show that $s_1<^{co}\gamma(s)$. Suppose not. Then $\gamma(s)<^{co} s_1\vee \gamma(s)$ and
\[ \vec{f}_{i-1}(s_1\vee \gamma(s))=\vec{f}_{i-1}(s_1)\vee \vec{f}_{i-1}(\gamma(s))\leq^{co} \gamma(t).\]
This means that
\[ \gamma(f^*_i(t))= \gamma(s) <^{co} s_1\vee \gamma(s) \leq^{co} f^*_{i-1}(\gamma(t)). \]
and that the codomains are not preserved. Thus $s_1<^{co}\gamma(s)$ indeed.

Next we show that $s_1\in \delta(s)$. Again, we suppose that this is not the case. Then there is $s_2^\circ\in \delta(s)$ such that $s_1<^{co}s^\circ_2$. We have
\[ \vec{f}_{i-1}(s_1)<^{co} \vec{f}_{i-1}(s_2^\circ) <^{co} \vec{f}_{i-1}(\gamma(s))\leq^{co} \gamma(t),   \]
and
\[ \vec{f}_{i-1}(s_1)\leq^{co} t_1 <^{co} \gamma(t). \]
As $f^*_{i-1}(t_1)=s_1$, we have
\[ \vec{f}_{i-1}(s_1)<^{co} t_1 <^{co}\vec{f}_{i-1}(s_2^\circ) <^{co} \vec{f}_{i-1}(\gamma(s))\leq^{co} \gamma(t). \]
As the set $\lvs^i(t)^\circ$ is a convex subtree of $T^*_{i-1}$, we have $\vec{f}_{i-1}(s_2^\circ)\in \lvs(t)^\circ$. Hence $s_2^\bullet\not\in\lvs^i(s)$ and $\vec{f}_i(s_2^\bullet)\in \lvs^i(t)$. Thus we have
\[ \vec{f}_{i}(s_2^\bullet\vee s)=\vec{f}_{i}(s_2^\bullet)\vee \vec{f}_i(s) \leq^{co} t, \]
and $s<^{co}s_2^\bullet\vee s$. This contradicts the fact $f^*_i(t)=s$. Thus $s_1 \in \delta(s)$, as claimed.

So far we have shown that $f^*_i$ restricts to a well defined function
\[ f^*_{i-1\lc t} :\delta(t)-\ker(f^*) \lra \delta(s). \]
We shall show that it is a bijection.

Let $t_1,t_2\in \delta(t)$ and $s\in S^*_{i-1}$ and $f^*_i(t_1)=s_1=f^*_i(t_2)$. Hence  $\vec{f}_{i-1}(s_1)<^{co}t_1$ and $\vec{f}_{i-1}(s_1)<^{co}t_2$ and then $t_1\perp^+t_2$ or $t_1=t_2$.
As $\delta(t)$ is an antichain in $S^*_{i-1}$, $t_1=t_2$. Thus $f^*_{\lc t}$ is one-to-one.

To see  that $f^*_{\lc t}$ is onto, let us fix an arbitrary $s_1\in \delta(s)$. Then $s_1\not\in\lvs(s)^\circ$.

We shall show that $\vec{f}_{i-1}(s_1)\not\in \lvs^i(t)^\circ$. Suppose to the contrary that $s_1=s_3^\circ$, for some $s_3\in S_{i}$ and that

\[ {\rm lvs}^i(t)^\circ \ni \vec{f}_{i-1} (s_3^\circ)\not\in \vec{f}_{i-1} (\lvs^i(s)^\circ).\]

Hence $s_3^\bullet\not\in\lvs^i(s)$ and $f(s_3^\bullet)\in \lvs^i(t)$. Thus

\[ \vec{f}_i(s_3^\bullet\vee s) =\vec{f}_i(s_3^\bullet)\vee \vec{f}_i(s)\leq^{co} t\]

and $s<s_3^\bullet\vee s$. This contradicts the fact that $f^*_i(t)=s$. Thus $\vec{f}_i(s_1)\not\in \lvs^i(t)^\circ$ indeed.

There is $t_1\in \delta(t)$ such that $\vec{f}_{i-1}(s_1)\leq^{co} t_1$. Let $s_2\in \lvs^i(s)^\circ$ such that $s_1\prec^{co} s_2$. Then
\[  \vec{f}_{i-1}(s_2) \in \vec{f}_{i-1}(\lvs^i(s)^\circ) \subseteq \lvs^i(t)^\circ.\]
Hence $s_1$ is the largest element of $S^*_{i-1}$ such that $\vec{f}_i(s_1)\leq^{co} t$, and hence $f^*_i(t_1)=s_1$, as required.
This ends the proof of Case $\delta.1$.

\vskip 2mm
{\em Case $\delta$.2}:  $f^*_i(t)=s_1\in  S^*_{i-1}$.

In this case we have a $t_1\in \delta(t)$ such that $\vec{f}_{i-1}(s_1)\leq^{co} t_1$.
Clearly $f^*_{i-1}(t_1)=s_1=f^*_{i-1}(t)$. It remains to show that $\delta(t)-\{t_1\}\subseteq \ker(f^*)$. Suppose to the contrary that there is $t_2\in \delta(t)$, $t_2\neq t_1$ such that $f^*_{i-1}(t_2)\in S^*_{i-1}$. Thus there is $s_2\in S^*_{i-1}$ such that $\vec{f}_{i-1}(s_2)\leq^{co}t_2\leq \gamma(t)$. Hence $s_1<s_1\vee s_2$ and
\[ \vec{f}_{i-1}(s_1\vee s_2) =  \vec{f}_{i-1}(s_1)\vee \vec{f}_{i-1}(s_2) \leq^{co} \gamma(t). \]
But then
\[ \gamma^{(i-1)}(f^*_i(t)) = \gamma^{(i-1)}(s_1)=s_1<^{co}s_1\vee s_2 \leq^{co} f^*_i (\gamma^{(i-1)}(t)) \]
and this contradicts the fact that $f^*$ preserves codomains. This ends the proof of Case $\delta.2$.

\vskip 2mm
{\em Case $\delta$.3}: $f^*_i(t)=s_1\in  S^*_{j}$, for some $j<i-1$.

We need to show that $\delta(t)\subseteq \ker(f^*)$. Suppose not. Then there is $t_1\in\delta(t)$ and $s_2\in S^*_{i-1}$ such that $\vec{f}_i(s_2)\leq^{co} t_1$. But this means
that $f^*_i(t)\in  S^*_{j}$, for some $j\geq i-1$, contrary to the supposition. $~~\Box$

For the proof of duality we need the following observations. We use the notation introduced above.

\begin{lemma}\label{lem-ope-obs}
Let $P$ be a positive opetope, $(P^*,\pi)$ corresponding positive zoom complex, $i\in\o$, $p\in P^*_{i+1}= P_{i+1}-\gamma(P_{i+2})$, $p_{root}= \sup_{<^-}(\pi(\gamma(p)))$. Then
\begin{enumerate}
  \item $\gamma\gamma(p)=\gamma(p_{root})$;
  \item\label{xyz} the map
  \[ \xi_p :(\pi(\gamma(p)),<^-) \lra (\lvs^i(p^\circ)^\circ,<^{co}),\] such that, for \mbox{$q\in \pi(\gamma(p))\subseteq P_{i}-\gamma(P_{i+1})$, $\xi_p (q)=q^\circ$ } is an order isomorphism.
  \item In particular, $p_{root}^\bullet =\xi_p (p_{root})=\sup^{i}(\lvs^{i+1}(p^\circ)).$
\end{enumerate}
\end{lemma}
{\it Proof.}~ Exercise. For \ref{xyz}. use Path Lemma. $~~\Box$

\newpage
\subsection{The main theorem}

In this section we shall prove that the functors defined in previous sections are essential inverse one to the other.

\begin{theorem}\label{duality}
The functors
\begin{center} \xext=1200 \yext=250
\begin{picture}(\xext,\yext)(\xoff,\yoff)
\putmorphism(0,100)(1,0)[\pOpe_{\iota,epi}`\pZoom^{op}`]{800}{0}a
\putmorphism(0,150)(1,0)[\phantom{\pOpe_{\iota,epi}}`\phantom{\pZoom^{op}}`(-)^*]{800}{1}a
\putmorphism(0,50)(1,0)[\phantom{\pOpe_{\iota,epi}}`\phantom{\pZoom^{op}}`(-)^*]{800}{-1}b
 \end{picture}
\end{center}
defined above, establish a dual equivalence of categories between categories of positive opetopes with $\iota$-epimorphisms and positive zoom complexes with embeddings.
\end{theorem}

{\it Proof.}~ We shall define two natural isomorphisms $\eta$ and $\varepsilon$.

Let $(S,\sigma)$ be a positive zoom complex.
Recall that \[ S^*_i=S_{i+1}\lhd S_i,\;\; {\rm and}\;\; S^{**}_i= (S_{i+1}\lhd S_i) - \gamma(S_{i+2}\lhd S_{i-1}).\]
For $i\in\o$, the $i$-th component
\[ \eta_{S,i}: S_i \lra S^{**}_i \]
of $\eta_S: (S,\sigma) \lra (S^{**},\sigma^{**})$ is defined as
\[ \eta_{S,i}(s) = s^\bullet,\]
i.e. it is a vertex in $S_{i+1}\lhd S_i$. Clearly $\eta_{S,i}$ is one-to-one and, as all circles in $S_{i+1}\lhd S_i$ are of form $\gamma(S_{i+2}\lhd S_{i-1})$, $\eta_{S,i}(s)$ is onto, as well. To see that  $\eta_{S,i}$ is an order isomorphism, consider $s_1,s_2\in S_i$. Then
\[ s^\bullet_1\prec^-s^\bullet_2 \]
iff
\[ \gamma(s^\bullet_1)\in \delta(s^\bullet_2) \]
iff
\[ s^\circ_1\in \cvr^{i-1}(s^\circ_2) \]
iff
\[ s^\circ_1\prec^{co} s^\circ_2 \]
iff
\[ s_1 \prec^{S_i} s_2. \]
To see that $\eta_{S}$ is an isomorphism of positive zoom complexes, it is enough to show that, for $i\in \o$,
$$(\eta_{S,i+1},\eta_{S,i}): (S_{i+1},\sigma_i,S_i)\lra (S^{**}_{i+1},\sigma^{**}_i,S^{**}_i)$$
is an isomorphism of constellations. To this aim, it is enough to show that the square
\begin{center} \xext=750 \yext=550
\begin{picture}(\xext,\yext)(\xoff,\yoff)
\setsqparms[1`1`1`1;700`500]
 \putsquare(0,0)[S_{i+1}`St(S_i)`S^{**}_{i+1}`St(S^{**}_{i});\sigma_i`\eta_{S,i+1}`\vec{\eta}_{S,i}`\sigma^{**}_i]
 \end{picture}
\end{center}
commutes, where $\vec{\eta}_{S,i}$ is the image function induced by the function $\eta_{S,i}$.  Let $s\in S^*_{i+1}$. We have

\[ \vec{\eta}_{S,i}(\sigma_i(s)) = \{ t^\bullet :\;  t\in S_i,\;  t \in \sigma_i(s) \} =   \]
\[ = \{ t^\bullet\in S^{**}_i : t^\bullet <^{co} s^\circ \} =   \]
\[ = \{ t^\bullet\in S^{**}_i : t^\bullet <^+ \gamma(s^\bullet) \} = \]
\[ = \sigma^{**}_i(s^\bullet) =  \sigma^{**}_i(\eta_{S,i}(s). \]
The naturality of $\eta$ is clear.

\vskip 3mm

Now we shall check that $\varepsilon$ is a natural isomorphism. Let $P$ be a positive opetope, $i\in \o$. By Lemma \ref{lem-eps-adjoint}, the maps
\[ \varepsilon_{P,i} : (P^*_{i+1}\lhd P^*_i,<^{co}) \lra (P_i,<^+) \]
defined in section \ref{from p-ope to t-compl} are order isomorphism. Recall that, for $p_1\in P_{i}-\gamma(P_{i+1})$, $p_1^\bullet \in P^*_{i+1}\lhd P^*_{i}$, we have $\varepsilon_{P,i}(p_1^\bullet)=p_1$ and, for $p_2\in P_{i+1}-\gamma(P_{i+2})$, $p_2^\circ \in P^*_{i+1}\lhd P^*_{i}$, we have $\varepsilon_{P,i}(p_2^\circ)=\gamma(p_2)$.

%
%

We need to show that $\varepsilon_{P}$ preserves both codomains $\gamma$ and domains $\delta$. Nautrality of $\varepsilon$ is again clear.

Preservation of codomains. Let $p_1\in P_{i+1}-\gamma(P_{i+2})$. We have
\[ \varepsilon_{P,i}(\gamma(p_1^\bullet)) = \varepsilon_{P,i}(p_1^\circ) =\]
\[   \gamma(p_1) = \gamma(\varepsilon_{P,i+1}(p_1^\bullet)). \]
Let $p_2\in P_{i+2}-\gamma(P_{i+3})$ and $p_{root}\in P_{i+1}-\gamma(P_{i+2})$ such that $p_{root}^\circ= \sup^i(\lvs^{i+1}(p_2^\circ)^\circ)$. Using Lemma \ref{lem-ope-obs},  we have
\[ \varepsilon_{P,i}(\gamma(p_2^\circ)) = \varepsilon_{P,i}(\sup^i(\lvs^{i+1}(p_2^\circ)^\circ)) = \]
\[ = \varepsilon_{P,i}(p^\circ_{root}) = \gamma(p_{root})=\]
\[ =\gamma\gamma(p_2) = \gamma(\varepsilon_{P,i+1}(p_2^\circ)). \]

Preservation of domains. Let $p\in P_{i+1}^{**}$ and $q\in P^{**}_i$. We need to verify that
\begin{equation}\label{eps-pres-doms}
   q\in\delta(p) \;\;\; {\rm iff}\;\;\; \varepsilon_{P,i}(q)\in\delta(\varepsilon_{P,i+1}(p)).
\end{equation}
We shall prove the above equivalence by cases depending on the form of $p$ and $q$.

Let $p_1\in P_{i+1}-\gamma(P_{i+2})$, $p_2\in P_{i+2}-\gamma(P_{i+3})$, $q_1\in P_{i}-\gamma(P_{i+1})$, $q_2\in P_{i+1}-\gamma(P_{i+2})$. Then we shall
consider four cases, one by one.

\vskip 2mm
{\em Case 1:} $p=p_1^\bullet$, $q=q_1^\bullet$. We have

\[ q_1^\bullet\in\delta(p_1^\bullet) \]
iff
\[ q_1^\bullet\in\cvr^{i-1}(p_1^\circ) \]
iff
\[ q_1^\bullet \prec^{co} p_1^\circ \]
iff
\[ q_1\in \delta(p_1) \]
iff
\[ \varepsilon_{P,i}(q_1^\bullet)\in\delta(\varepsilon_{P,i+1}(p_1^\bullet)). \]

\vskip 2mm
{\em Case 2:} $p=p_1^\bullet$, $q=q_2^\circ$.

\[ q_2^\circ\in\delta(p_1^\bullet) \]
iff
\[ q_2^\circ\prec^{co} p_1^\circ \]
iff
\[ q_2\prec^{-} p_1 \]
iff
\[ \gamma(q_2)\in \delta( p_1) \]
iff
\[ \varepsilon_{P,i}(q_2^\circ)\in\delta(\varepsilon_{P,i+1}(p_1^\bullet)). \]

\vskip 2mm
{\em Case 3:} $p=p_2^\circ$, $q=q_1^\bullet$.

\[ q_1^\bullet\in\delta(p_2^\circ) \]
iff
\[ q_1^\bullet\in\cvr^{i-1}(\lvs^1(p_2^\circ)^\circ) \]
iff
\[ q_1^\bullet\not\in\lvs^1(p_2^\circ)^\circ\;{\rm and}\; {\rm there\; is}\; q_3\in P_i-\gamma(P_{i+1})\; {\rm such\; that}\; q_1^\bullet\prec^{co} q_3^\circ\; {\rm and}\;q_3^\circ \in \lvs^i(p_2^\circ)^\circ  \]
iff
\[ {\rm there\; is}\; q_3\in P_i-\gamma(P_{i+1})\; {\rm such\; that}\; q_1\in \delta(q_3) \; {\rm and}\; q_3\leq^+ \gamma(p_2)  \]
iff (Path Lemma)
\[ q_1\in \delta\gamma(p_2) \]
iff
\[ \varepsilon_{P,i}(q_1^\bullet)\in\delta(\varepsilon_{P,i+1}(p_2^\circ)). \]

\vskip 2mm
{\em Case 4:} $p=p_2^\circ$, $q=q_2^\circ$.

\[ q_2^\circ\in\delta(p_2^\circ) \]
iff
\[ q_2^\circ\in\cvr^{i-1}(\lvs^1(p_2^\circ)^\circ) \]
iff
\[ q_2^\circ\not\in\lvs^1(p_2^\circ)^\circ\;{\rm and}\; {\rm there\; is}\; q_3\in P_{i+1}-\gamma(P_{i+2})\; {\rm such\; that}\; q_2^\circ\prec^{co} q_3^\circ\; {\rm and}\; q_3^\circ \in \lvs^i(p_2^\circ)^\circ  \]
iff
\[ q_2\not\leq^+\gamma(p_2) \;{\rm and}\; {\rm there\; is}\; q_3\in P_{i+1}-\gamma(P_{i+2})\; {\rm such\; that}\; q_2\prec^{-} q_3\; {\rm and}\; q_3 \leq^+\gamma(p_2)  \]
iff (Path Lemma)
\[ \gamma(q_2)\in \delta\gamma(p_2) \]
iff
\[ \varepsilon_{P,i}(q_2^\circ)\in\delta(\varepsilon_{P,i+1}(p_2^\circ)). \]
$~~\Box$


\section{Wide zoom complexes and opetopic cardinals}\label{sec-wide-complex}

In this section we extend the above duality to the wide zoom complexes with embeddings on one side and positive opetopic cardinals with $\iota$-epimorphisms on the other.

\subsection{Wide constellations}

A {\em forest} is a finite poset $(P,\leq)$ which is a disjoint sum of trees. A morphism of forests $f:(P,\leq)\ra (Q,\leq)$  is a one-to-one function that
preserves and reflects the order. $\St(P)$ is the poset of convex sub-trees of the forest $P$.

A {\em wide constellation }
is a triple $(T',\tau,T)$ where $T$, $T'$ are forests and $\tau$ is a monotone function
\[ \tau: T'\ra \St(T) \]
such that if $t,t'\in T'$ and $\sigma(t)\cap \sigma(t')\neq\emptyset$, then $t \perp t'$.

Let $\tau: T'\ra \St(T)$ be a constellation. Then the {\em constellation forest} $T'\lhd_\tau T$, is the extension of $T'$ by $T$ along $\tau$, i.e., it is the forest $T'$ with nodes of $T$ added as leaves so that if $x\in T$ and $y\in T'$, then $x<^{co}y$ in $T'\lhd_\tau T$ iff $x\in \tau(y)$. The order $<^{co}$ is called the {\em constellation order of the wide constellation} $\tau: T'\ra \St(T)$, or just the {\em constellation order} if the constellation is understood. Any pair of maps of forests $f:S\ra T$, $f':S'\ra T'$ such that $\vec{f}(\sigma(s))\subseteq \tau(f'(s))$ for $s\in S$, induces a monotone map $f' \lhd f: S'\lhd_\sigma S\lra T'\lhd_\tau T$. Such a pair
\[ (f',f) : (S',\sigma,S) \lra (T',\tau,T)\]
is {\em a morphism of constellations} iff the induced map of constellation forests  $f' \lhd f$ preserves (existing) binary sups.

\subsection{Wide zoom complexes and duality}
{\em A wide zoom complex} $(T,\tau)$ is a sequence of wide constellations:
\[ \tau_0: T_1 \ra \St(T_{0}), \]
\[ \tau_1: T_2 \ra \St(T_1), \]
\[ \ldots \]
\[ \tau_{i}: T_{i+1} \ra \St(T_i), \]
\[ \ldots \]
$i\in\o$, with almost all sets $T_i$ empty. The dimension $(T,\tau)$ is $n$ iff  $T_n$ is the last non-empty set. We write $dim(T)$ for dimension of the wide zoom complex $(T,\tau)$. $T_0$ is required to be a singleton.

A {\em morphism of wide zoom complexes} $f:(S,\sigma)\ra (T,\tau)$ is a family of forest embeddings $f_i:S_i\ra T_i$, for $i\in \o$, such that,
\[(f_{i+1},f_i): (S_{i+1},\sigma_{i},S_i)\lra (T_{i+1},\tau_{i},T_i)\]
is a morphism of constellations, for $i\in \o$.

The {\em size of a wide zoom complex} $(T,\tau)$ is a sequence of natural numbers $size(T,\tau)=\{s_i\}_{i\in \o}$ so that $s_i=size(T,\tau)_i$ is the number of trees in the forest $T_i$. A wide zoom complex $(T,\tau)$ is a positive zoom complex iff $size(T,\tau)_i\leq 1$, for all $i\in\o$.

The category of wide zoom complexes and their morphisms will be denoted by $\wZoom$. Clearly $\pZoom$ is a full subcategory of $\wZoom$.
\begin{theorem}
The functors
\begin{center} \xext=1200 \yext=250
\begin{picture}(\xext,\yext)(\xoff,\yoff)
\putmorphism(0,100)(1,0)[\pOpeCard_{\iota,epi}`\wZoom^{op}`]{1200}{0}a
\putmorphism(0,150)(1,0)[\phantom{\pOpeCard_{\iota,epi}}`\phantom{\wZoom^{op}}`(-)^*]{1200}{1}a
\putmorphism(0,50)(1,0)[\phantom{\pOpeCard_{\iota,epi}}`\phantom{\wZoom^{op}}`(-)^*]{1200}{-1}b
 \end{picture}
\end{center}
defined as those for trees, establish a dual equivalence of categories between categories of opetopic cardinals with $\iota$-epimorphisms and wide zoom complexes with embeddings.
\end{theorem}

{\it Proof}~ This is an easy extension of the corresponding fact concerning positive zoom complexes and positive opetopes.
$~~\Box$

\section{Appendix: Positive opetopes and positive opetopic cardinals}\label{sec-p-ope}

In this appendix we recall the notion of positive opetopes, positive opetopic cardinals, their morphisms: face maps and $\iota$-maps. We also quote without proofs some facts from
[Z1] and [Z4].

\subsection{Positive hypergraphs}

A {\em positive hypergraph} $S$ is a family $\{ S_k\}_{k\in \o}$ of finite sets of faces, a family of functions $\{ \gamma_k : S_{k+1}\ra S_k \}_{k \in\o }$, and a family of total relations $\{ \delta_k : S_{k+1}\ra S_k\}_{k \in \o}$. Moreover, $\delta_0 : S_1\ra S_0$ is a function and only finitely many among sets $\{ S_k\}_{k\in \o}$ are non-empty. As it is always clear from the context, we shall never use the indices of the functions $\gamma$ and $\delta$. We shall ignore the difference between $\gamma(x)$ and $\{ \gamma(x)\}$ and in consequence we shall consider iterated applications of $\gamma$'s and $\delta$'s as sets of faces, e.g. $\delta\delta(x)=\bigcup_{y\in \delta(x)} \delta(y)$ and  $\gamma\delta(x)=\{\gamma(y)\, |\, y\in \delta(x)\}$.

A {\em morphism of positive hypergraphs} $f:S\lra T$ is a family of functions $f_k : S_k \lra T_k$, for $k \in\o$, such that, for $k>0$ and $a\in S_k$, we have  $\gamma (f(a))=f(\gamma(a))$ and $f_{k-1}$ restricts to a bijection
\[ f_a: \delta(a) \lra \delta(f(a)). \]
The category of positive hypergraphs is denoted by $\bf pHg$.

We define a binary relation of {\em lower order} on $<^{S_k,-}$ for $k>0$ as the transitive closure of  the relation $\lhd^{S_k,-}$ on $S_k$ such that, for $a,b\in S_k$,  $a \lhd^{S_k,-} b$ iff $\gamma (a)\in \delta(b)$.  We write $a\perp^- b$ iff either $a <^- b$ or $b <^- a$, and we write $a \leq^- b$ iff either $a=b$ or $a <^- b$.

We also define a  binary relation of {\em upper order} on $<^{S_k,+}$ for $k\geq0$ as the transitive closure of  the relation $\lhd^{S_k,+}$ on $S_k$ such that, for $a,b\in S_k$,  $a \lhd^{S_k,+} b$ iff there is $\alpha\in S_{k+1}$ so that $a\in \delta(\alpha)$ and $\gamma (\alpha)=b$.  We write $a\perp^+ b$ iff either $a <^+ b$ or $b <^+ a$, and we write $a \leq^+ b$ iff either $a=b$ or $a <^+ b$.

\subsection{Positive opetopic cardinals}

A positive hypergraph $S$ is a {\em positive opetopic cardinal} if it is non-empty,
i.e. $S_0\neq\emptyset$ and it satisfies the following four conditions.
\begin{enumerate}
 \item {\em Globularity:}  for  $a\in S_{\geq 2}$:
 \[ \gamma\gamma(a) =\gamma\delta(a)-\delta\delta(a),\hskip 15mm \delta\gamma(a) =\delta\delta(a)-\gamma\delta(a).\]
  \item {\em Strictness:} for $k\in\o$, the relation   $<^{S_k,+}$ is a strict order; $<^{S_0,+}$ is linear.
  \item {\em Disjointness:} for $k>0$,
\[\perp^{S_k,-}\cap \perp^{S_k,+}=\emptyset.\]
  \item {\em Pencil linearity:} for any $k>0$   and $x\in S_{k-1}$, the sets
  \[ \{ a\in S_k \; | \; x=\gamma(a) \} \;\;\;\;{\rm   and}\;\;\;\; \{ a\in S_k \; | \; x\in \delta(a) \} \]
  are linearly ordered by $<^{S_k,+}$.
\end{enumerate}

The category of positive opetopic cardinals is the full subcategory of $\bf pHg$ whose objects are positive opetopic cardinals. It is denoted by $\pOpeCard$.

\subsection{Positive opetopes}

The {\em size of positive opetopic cardinal} $S$ is a sequence of natural numbers $size(S)=\{ | S_n - \delta (S_{n+1})| \}_{n\in\o}$, with
all elements above $dim(S)$ being equal $0$. We have an order $<$ on such sequences of natural numbers so that $\{ x_n \}_{n\in\o} < \{ y_n \}_{n\in\o}$ iff there is $k\in\o$ such that $x_k< y_k$ and, for all $l>k$, $x_l = y_l$. This order is well founded and hence facts about positive opetopic cardinals can be proven by induction on their size.

Let $P$ be an positive opetopic cardinal.  We say that $P$ is a {\em positive opetope} iff $size(P)_l\leq 1$, for $l\in \o$. By $\pOpe$ we denote full
subcategory of $\bf pHg$ whose objects are positive opetopes.

\vskip 2mm
{\bf Some notions and notation.}
\begin{enumerate}
\item Let $S$ be a positive hypergraph $S$, $x$ a face of $S$. 
The dimension of $S$ is maximal $k$ such that $S_k\neq\emptyset$. We denote by $dim(S)$ the dimension of $S$. We usually tacitly assume that the sets of faces of different dimensions are disjoint and we denote by $|S|=\bigcup_{i\in \o}S_i$ the sum of all faces of $S$.
  \item Let $P$ be an opetope.  If  $dim(P)=n$, then the unique face in $P_n$ is denoted by $\bm_P$.
  \item The function $\gamma^{(k)}: P\ra P_{\leq k}$ is an iterated version of the codomain function $\gamma$ defined as follows. For  any $k,l\in \o$ and $p\in P_l$,
\[ \gamma^{(k)}(p) \;= \; \left\{ \begin{array}{ll}
		\gamma\gamma^{(k+1)}(p)   & \mbox{ if } l>k \\
		p  & \mbox{ if } l\leq k.
                                    \end{array}
			    \right. \]
  \item   Let $a,b\in P_k$. A {\em lower path} $a_0, \ldots , a_m$ {\em from}\index{path!lower} $a$ {\em to} $b$ in $P$ is a sequence of faces $a_0, \ldots , a_m\in S_k$ such that $a=a_0$, $b=a_m$ and, for $\gamma(a_{i-1})\in\delta(a_i)$, $i=1,\ldots ,m$.
\item   Let $x,y\in P_k$. An  {\em upper path} $x,a_0, \ldots , a_m,y$ {\em from}\index{path!upper} $x$ {\em to} $y$ in $S$ is a  sequence of faces $a_0, \ldots , a_m\in P_{k+1}$ such that  $x\in\delta(a_0)$, $y=\gamma(a_m)$ and $\gamma(a_{i-1})\in\delta(a_i)$, for $i=1,\ldots ,m$.
\end{enumerate}

\begin{lemma}
\label{tech_lemma1} Let $P$ be a positive opetopic cardinal,
$n\in\o$, $a,b\in P_n$, $a<^+b$. Then, there is an upper
$P_{n+1}-\gamma(P_{n+2})$-path from $a$ to $b$. $\Box$
\end{lemma}

\begin{lemma}[Path Lemma]\label{fullpath} Let $P$ be an opetope. Let $k\geq 0$, $B=(a_0,\ldots, a_k )$ be a
maximal $S_n$-lower path in a positive opetopic cardinal $P$,  $b\in S_n$,
$0\leq s\leq k$, $a_s<^+b$. Then there are  $0\leq l\leq s\leq p
\leq k$ such that
\begin{enumerate}
   \item $a_i<^+b$ for $i= l, \ldots , p$;
   \item $\gamma(a_p)=\gamma(b)$;
  \item either $l=0$ and $\delta(a_0)\subseteq \delta(b)$ or
  $l>0$ and $\gamma(a_{l-1})\in\delta(b)$;
  \item $\gamma(a_i)\in\iota(S)$, for $l\leq i <p$. $\Box$
\end{enumerate}
\end{lemma}

\subsection{The embedding of \texorpdfstring{$\pOpeCard$}{pOpeCard} into \texorpdfstring{$\oCat$}{oCat}}

There is an embedding
\[ (-)^*:\pOpeCard \lra \oCat \]
defined as follows, c.f. \cite{Z1}. Let $T$ be an opetopic cardinal. The $\o$-category $T^*$ has as $n$-cells pairs $(S,n)$, where $S$ is a subopetopic cardinal of $T$, $dim(S)\leq n$, and $n\geq 0$.

For $k<n$, the domain and codomain operations
\[ \bd^{(k)}, \bc^{(k)}: T^*_{n} \lra T^*_k \]
are given, for  $(S,n)\in T^*_{n}$, by
\[  (\bd^{(k)}(S,n))=(\bd^{(k)}(S),k),\hskip 15mm  (\bc^{(k)}(S,n))=(\bc^{(k)}(S),k) \]
where
\[ (\bd^{(k)}(S)_l \;= \; \left\{ \begin{array}{ll}
		\emptyset   & \mbox{ if } l>k \\
		S_k-\gamma(S_{k+1})  & \mbox{ if } l=k \\
        S_l  & \mbox{ if } l<k
                                    \end{array}
			    \right. \]
and
\[ (\bc^{(k)}(S)_l \;= \; \left\{ \begin{array}{ll}
		\emptyset   & \mbox{ if } l>k \\
		S_k-\delta(S_{k+1})  & \mbox{ if } l=k \\
        S_{k-1}-\iota(S_{k+1})  & \mbox{ if } l=k-1\geq 0 \\
        S_l  & \mbox{ if } l<k-1
                                    \end{array}
			    \right. \]
The identity operation
\[ id : T^*_n \lra T^*_{n+1} \]
is given by
\[  (S,n)\mapsto (S,n+1). \]
The composition operation is defined, for pairs of cells $(S,n),(S',n')\in T^*$ with $k\leq n,n'$ such that $\bd^{(k)}(S,n)=\bc^{(k)}(S',n')$, as the sum
\[ (S,n)\circ (S',n')= (S\cup S', \max(n,n')).\]

Now $T^*$ together with operations of domain, codomain, identity and composition is an $\o$-category. If $f:S\ra T$ is a map of opetopic cardinals and $S'$ is a sub-opetopic cardinal of $S$, then the image $f(S')$ is a sub-opetopic cardinal of $T$. This defines the functor $(-)^*$ on morphisms. We recall from \cite{Z1}

\begin{theorem}\label{thm-embedding-pOpeCard}
The embedding
\[ (-)^*:\pOpeCard \lra \oCat \]
is well defined and full on isomorphisms and it factorises through $\Poly\lra \oCat$ via full and faithful functor, $(-)^*:\pOpeCard \lra \Poly$, into the category of polygraphs. $\Box$
 \end{theorem}

\subsection{\texorpdfstring{$\iota$}{i}-maps of positive opetopes}\label{sec-iota-maps}

The embedding $(-)^\ast: \pOpe \ra \oCat$ is not full, but it is full on isomorphisms. The morphisms $P^\ast\ra Q^\ast$ in $\oCat$ between images of opetopes are $\o$-functors that send generators to generators.  The category $\pOpei$ with the same objects as $\pOpe$ will be so defined that the embedding $(-)^\ast: \pOpei \ra \oCat$ (denoted the same way) will be full on $\o$-functors that send generators to either generators or identities on generators of a smaller dimension.

Let $P$ and $Q$ be positive opetopes. A {\em contraction morphism of opetopes} (or {\em $\iota$-map}, for short),  $h : Q \ra P$, is a function $h : |Q| \ra |P|$ between faces of opetopes such that
\begin{enumerate}
  \item  $dim(q)\geq dim(h(q))$, for $q\in  Q$;
  \item (preservation of codomains) $h(\gamma^{(k)} (q))= \gamma^{(k)}(h(q))$, for $k\geq 0$ and $q\in Q_{k+1}$;
  \item (preservation of domains)
  \begin{enumerate}
    \item  if $dim(h(q))= dim(q)$, then $h$ restricts to a bijection
  \begin{center}
\xext=1200 \yext=10 \adjust[`I;I`;I`;`I]
\begin{picture}(\xext,\yext)(\xoff,\yoff)
\putmorphism(0,0)(1,0)[(\delta(q)-\ker(h))`\delta(h(q))`h]{800}{1}a
\end{picture}
\end{center}
for $k\geq 0$ and $q\in Q_{k+1}$, where {\em the kernel of h} is defined as $$\ker(h)=\{ q\in Q| dim(q)> dim(h(q))\};$$
    \item if $dim(h(q))= dim(q)-1$, then $h$ restricts to a bijection
  \begin{center}
\xext=1200 \yext=10 \adjust[`I;I`;I`;`I]
\begin{picture}(\xext,\yext)(\xoff,\yoff)
\putmorphism(0,0)(1,0)[(\delta(q)-\ker(h))`\{ h(q)\}`h]{800}{1}a
\end{picture}
\end{center}
for $k\geq 0$ and $q\in Q_{k+1}$;
    \item if $dim(h(q))< dim(q)-1$, then $\delta^{(k)}(q)\subseteq \ker(h)$.
  \end{enumerate}
\end{enumerate}

We have an embedding $\kappa:\pOpe\lra \pOpei$ that induces the usual adjunction $\kappa_!\dashv \kappa^*$

\begin{center}
\xext=1000 \yext=370
\begin{picture}(\xext,\yext)(\xoff,\yoff)
\putmorphism(0,140)(1,0)[\widehat{\pOpe}`\widehat{\pOpei}`]{1000}{0}a
\putmorphism(0,210)(1,0)[\phantom{\widehat{\pOpe}}`\phantom{\widehat{\pOpei}}`\kappa_!]{1000}{1}a
\putmorphism(0,70)(1,0)[\phantom{\widehat{\pOpe}}`\phantom{\widehat{\pOpei}}`\kappa^*]{1000}{-1}b
\end{picture}
\end{center}

\begin{lemma}Let $h: Q\ra P$ be a $\iota$-map, $q_1,q_2\in Q-\ker(h)$ and $l<k\in \o$ such that
\[  \begin{array}{rcl}
		\gamma^{(k+1)}(q_1)  & <^- & \gamma^{(k+1)}(q_2) \\
		\gamma^{(k)}(q_1)   & <^+ & \gamma^{(k)}(q_2)\\
        \ldots &\ldots & \ldots \\
        \gamma^{(l+1)}(q_1)   & <^+ & \gamma^{(l+1)}(q_2)\\
        \gamma^{(l)}(q_1)   & = & \gamma^{(l)}(q_2).\\
                                    \end{array}
			    \]
Then there is $l\leq l'< k$ such that
\[  \begin{array}{rcl}
		h(\gamma^{(k+1)}(q_1))  & <^- & h(\gamma^{(k+1)}(q_2)) \\
		h(\gamma^{(k)}(q_1))   & <^+ & h(\gamma^{(k)}(q_2))\\
        \ldots &\ldots & \ldots \\
        h(\gamma^{(l'+1)}(q_1))   & <^+ & h(\gamma^{(l'+1)}(q_2)) \\
        h(\gamma^{(l')}(q_1))   & = & h(\gamma^{(l')}(q_2)).\; \Box\\
                                    \end{array}
			    \]
\end{lemma}

From the above we get immediately

\begin{corollary}\label{coro-iota-preservation}
Let $h: Q\ra P$ be a $\iota$-map, $q_1,q_2\in Q-\ker(h)$. Then
\begin{enumerate}
  \item $q_1<^- q_2$ iff $h(q_1)<^- h(q_2)$;
  \item if $q_1<^+ q_2$, then $h(q_1)\leq^+ h(q_2)$;
  \item if $h(q_1)<^+ h(q_2)$, then $q_1<^+ q_2$;
  \item if $h(q_1)= h(q_2)$, then $q_1\perp^+ q_2$. $\Box$
\end{enumerate}
\end{corollary}

A set $X$ of $k$-faces in a positive opetope $P$ is a {\em $<^+$-interval} (or {\em interval}, for short)  if it is either empty or there are two k-faces $x_0,x_1\in P_k$ such that $x_0\leq^+ x_1$ and $X=\{ x\in P_k | x_0\leq^+ x \leq^+ x_1 \}$. Any interval in any positive opetope is linearly ordered by $\leq^+$.

\begin{corollary}
Let $h:Q\ra P$ -be a contraction of positive opetopes, $p\in P_k$. Then the fiber of $k$-faces (of non-degenerating faces) $h^{-1}(p)- \ker(h)$ of $h$ over $p$ is an interval.
\end{corollary}

{\em Proof.} From Corollary \ref{coro-iota-preservation}.4 we get that $h^{-1}(p)- \ker(h)$ is linearly ordered. And from Corollary \ref{coro-iota-preservation}.2 that this linear order is an interval. $\Box$

\subsection{The embedding of $\pOpei$ into $\oCat$}

We extend the embedding functor $(-)^\ast$ to contractions
\[ (-)^* : \pOpei \lra \oCat. \]
Let $h:Q \ra P$ be a contraction morphism in $ \pOpei$. Then
\[ h^\ast : Q^\ast\ra P^\ast \]
is an $\o$-functor such that
\[ h^\ast(k,A) =  (k,\vec{h}(A)) \]
where $(k,A)\in Q^*_k$, and $\vec{h}(A)$ is the set-theoretic image of the opetopic cardinal $A$ under $h$.
\vskip 2mm

\begin{theorem}\label{theorem-iota-emdedding}The functor
\[ (-)^*: \pOpei \lra \oCat\]
is well defined. The objects of $\pOpei$ are sent under $(-)^*$ to positive-to-one polygraphs.  $(-)^*$ is faithful, conservative and full on those $\o$-functors that send generators to either generators or to (possibly iterated) identities on generators of smaller dimensions. In particular, it is full on isomorphisms. $\Box$
\end{theorem}

\end{document}